\documentclass[12pt, reqno]{amsart}
\usepackage{amsmath}
\usepackage{amssymb}
\usepackage{epsfig}
\usepackage{graphicx}
\usepackage{color}
\usepackage{fullpage}
\definecolor{shadecolor}{gray}{0.875}
\usepackage{amscd}
\usepackage{comment}
\usepackage{enumitem}
\usepackage{mathrsfs}

\usepackage[colorlinks=true]{hyperref}

\numberwithin{equation}{section}

\input xy
\xyoption{all}

\calclayout
\allowdisplaybreaks[3]

\theoremstyle{plain}
\newtheorem{prop}{Proposition}[section]

\newtheorem{theo}[prop]{Theorem}
\newtheorem{coro}[prop]{Corollary}

\newtheorem{lemm}[prop]{Lemma}

\theoremstyle{definition}
\newtheorem{defi}[prop]{Definition}

\newtheorem{conj}[prop]{Conjecture}

\newtheorem{rema}[prop]{Remark}

\newtheorem{exam}[prop]{Example}

\newtheorem{nota}[prop]{Notation}

\def\lra{\longrightarrow}

\def\cC{{\mathcal C}}

\def\cM{{\mathcal M}}
\def\cN{{\mathcal N}}
\def\cO{{\mathcal O}}

\def\cT{{\mathcal T}}

\def\cX{{\mathcal X}}

\def\ocM{\overline{\mathcal M}}

\def\bR{{\mathbb R}}
\def\bZ{{\mathbb Z}}

\def\Pic{\mathrm{Pic}}

\def\Supp{\mathrm{Supp}}

\def\Spec{\mathrm{Spec}}

\def\Pic{\mathrm{Pic}}

\def\Sec{\mathrm{Sec}}

\def\ch{\operatorname{char}}
\def\cok{\operatorname{Cok}}
\def\ev{\operatorname{ev}}

\makeatother
\makeatletter

 \author{Qile Chen}
\address{Department of Mathematics \\
Boston College  \\
Chestnut Hill, MA \, \, 02467}
\email{qile.chen@bc.edu}

\author{Brian Lehmann}
\address{Department of Mathematics \\
Boston College  \\
Chestnut Hill, MA \, \, 02467}
\email{lehmannb@bc.edu}

\author{Sho Tanimoto}
\address{Graduate School of Mathematics, Nagoya University, Furocho Chikusa-ku, Nagoya, 464-8602, Japan}
\email{sho.tanimoto@math.nagoya-u.ac.jp}

\title[Campana rational connectedness]{Campana rational connectedness\\ and weak approximation}

\begin{document}
\date{\today}

\begin{abstract}
Campana introduced a notion of Campana rational connectedness for Campana orbifolds. 
Given a Campana fibration over a complex curve, we prove that a version of weak approximation for Campana sections holds at places of good reduction when the general fiber satisfies a slightly stronger version of Campana rational connectedness.
Campana also conjectured that any Fano orbifold is Campana rationally connected; we verify a stronger statement for toric Campana orbifolds. A key tool in our study is log geometry and moduli stacks of stable log maps. 
\end{abstract}

\maketitle

\setcounter{tocdepth}{1}
\tableofcontents

\section{Introduction}

One of the major goals in arithmetic algebraic geometry is to understand rational points on algebraic varieties defined over the function field of a smooth projective curve over an algebraically closed field. By the valuative criterion, this amounts to studying the spaces of sections of fibrations over curves. 
In characteristic $0$ \cite[Theorem 1.1]{GHS03} showed that any rationally connected fibration admits a section,
and moreover by \cite[Theorem 3]{HT06} we know that weak approximation holds for such fibrations at places of good reduction.

On the other hand, in arithmetic geometry a major topic of interest is semi-integral points on Campana orbifolds.  Recent activity has focused on the notion of Campana points which interpolate between rational points and integral points. Many conjectures on rational/integral points admit a version for Campana points, e.g., \cite{AVA18, PSTVA, MNS24}.  Some recent works analyzing these conjectures include \cite{BY21, Xiao, Streeter22, Shute22, PS24, NS24, CLTBT24, Moerman}. However, there is still only limited evidence supporting these conjectures and more investigation is required to formulate them precisely. In this paper, we study Campana curves/sections and explore the validity of weak approximation in the setting of Campana sections. 

We first lay down the foundation of Campana curves/sections using log geometry and moduli stacks of stable log maps. Then we reinterpret the notion of orbifold rational connectedness introduced by Campana and show that this notion is equivalent to a stronger property called Campana rational connectedness.
Then we prove that for any fibration whose general fibers satisfy a version of Campana rational connectedness, weak approximation for Campana sections at places of good reduction holds. Finally we verify this property for toric Campana orbifolds.

\subsection{Campana curves}

Throughout the paper, our ground field is an algebraically closed field $\mathbf k$ of arbitrary characteristic $p= \mathrm{char} \, \mathbf k$, but for simplicity we assume that $\mathbf k$ has characteristic $0$ in this introduction.

Let $\underline{X}$ be a smooth projective variety equipped with a strict normal crossings (SNC) divisor $\Delta = \cup_i \Delta_i$. Let $X$ be the log scheme associated to $(\underline{X}, \Delta)$. (See Section~\ref{ss:log-map-stacks} for its definition.) To each irreducible component $\Delta_{i}$ we assign a weight
\[
\epsilon_i = 1- \frac{1}{m_i},
\]
where $m_i \geq 1$ is an integer.  We then define the effective $\mathbb Q$-divisor
\[
\Delta_\epsilon = \sum_i \epsilon_i \Delta_i.
\]
The pair $(X, \Delta_\epsilon)$ is called a klt Campana orbifold (in the sense of \cite[Section 1.2.1]{Campana04}).
Note that for a Campana orbifold we always assume the pair $(X,\Delta)$ is log smooth; we include the term klt to emphasize that the coefficients of $\Delta_{\epsilon}$ are smaller than $1$.

\begin{defi}
Let $(\pi : C \to S, f : C \to X)$ be a stable log map with the canonical log structure such that $S$ is a geometric log point with the trivial log structure.  (This implies that the underlying scheme $\underline{C}$ is irreducible and smooth and the image $f(C)$ is not contained in the boundary $\Delta$.) Let $p_k$ be a marked point and $\mathbf{c}_k = (c_{k, i})$ be the contact order, i.e., $c_{k, i}$ is the local multiplicity of $f^*\Delta_i$ at $p_k$. We say $f : C \to X$ is a Campana curve for $(X, \Delta_\epsilon)$ if $c_{k, i} \geq m_i$ whenever $c_{k, i} \neq 0$.  
\end{defi}

Log geometry controls deformations of stable log maps.  Since log deformations keep the contact orders constant, log geometry is particularly useful to understand deformations of Campana curves. In this paper, we lay down the foundation of Campana curves using log geometry.

A central theme in our work is the existence of Campana rational curves on Campana orbifolds.  The following definition is a variant of the notions of orbifold uniruledness and orbifold rational connectedness pioneered by Campana, e.g., \cite{Campana07, Campana10, Campana11}:

\begin{defi} \label{defi:introrc}
Let $(X, \Delta_\epsilon)$ be a klt Campana orbifold. We say $(X, \Delta_\epsilon)$ is Campana uniruled if there is a dominant family of genus $0$ Campana curves whose underlying curves have non-trivial numerical class. Moreover, if any two general points on $X$ are contained in a family of Campana curves of genus $0$ that dominates $X$, we say $(X, \Delta_\epsilon)$ is Campana rationally connected.
\end{defi}

As in the case of rational curves, we prove that these two notions are respectively equivalent to the existence of free or very free Campana rational curves in characteristic $0$.

\begin{rema}
The ``orbifold'' notions due to Campana are presented differently but turn out to be equivalent to Definition \ref{defi:introrc}.  We discuss the relationship in Remark \ref{rema:campanasdef}.
\end{rema}

In \cite[Section 5.4]{Campana07}, Campana made several precise conjectures about the relationship between Campana uniruledness and the behavior of the orbifold tangent bundle.  We will primarily be interested in the following conjecture, which is a special case of \cite[Conjecture 9.10]{Campana11} and provides the main source of examples of Campana rational connectedness.

\begin{conj}
\label{conj:main_intro}
Assume that $\mathbf k$ has characteristic $0$.
Let $(X, \Delta_\epsilon)$ be a klt Fano orbifold, i.e, $(X, \Delta_\epsilon)$ is a klt Campana orbifold and $-(K_{\underline{X}} + \Delta_\epsilon)$ is ample. Then $(X, \Delta_\epsilon)$ is Campana rationally connected. 
\end{conj}

\begin{exam}[Campana]
Here we introduce examples of Campana rationally connected orbifolds found by Campana.
Let $(X, \Delta_\epsilon)$ be a klt Fano orbifold such that all irreducible components $\Delta_{i}$ have the same multiplicity $m \in \mathbb{Z}_{\geq 1}$.
Assume that the boundary divisor $\Delta = \sum_i \Delta_i$ is divisible by $m$ in $\Pic(\underline{X})$. Let $\rho : \underline{Y} \to \underline{X}$ be the degree $m$ cyclic cover of $X$ totally ramified along $\Delta$.
Then since we have
\[
-K_{\underline{Y}} = -\rho^*(K_{\underline{X}} + \Delta_\epsilon),
\]
the variety $\underline{Y}$ is a klt Fano variety. By \cite[Theorem 1]{Zhang06} or \cite[Corollary 1.3 and Corollary 1.5]{HM07} 
$\underline{Y}$ is rationally connected.
Then note that any rational curve on $\underline{Y}$ is a Campana rational curve on $(X, \Delta_\epsilon)$ after imposing the canonical log structure, so in particular $(X, \Delta_\epsilon)$ is Campana rationally connected. % , and moreover when $X$ has Picard rank $1$, $X$ is strongly Campana uniruled.
Such a construction applies when $\underline{X}$ is a smooth Fano complete intersection in $\mathbb P^n$ and $\Delta$ is a SNC divisor which is the restriction of a Cartier divisor on $\mathbb P^n$.
\end{exam}

\subsection{Main results} Let $B$ be a smooth projective curve defined over $\mathbf k$ with the trivial log structure.
Let $\underline{\mathcal X}$ be a smooth projective variety equipped with a flat morphism $\pi : \underline{\mathcal X} \to B$ whose fibers are connected.
Let $\Delta = \cup_i \Delta_i$ be a SNC divisor on $\underline{\mathcal X}$ such that $\pi|_{\Delta} : \Delta \to B$ is flat and let $\mathcal X$ be the log scheme associated to $(\underline{\mathcal X}, \Delta)$.
By a klt Campana fibration $(\mathcal X/B, \Delta_\epsilon)$, we mean the data of a klt Campana orbifold $(\mathcal X, \Delta_\epsilon)$ equipped with a fibration $\pi: \underline{\mathcal X} \to B$ as described above.  
In the setting of a klt Campana fibration, one is interested in the moduli space of log sections $\sigma: C \to \mathcal{X}$ which satisfy the Campana condition.

A Campana jet for $(\mathcal{X},\Delta_{\epsilon})$ is a jet whose local multiplicities satisfy the Campana condition.  (See Definition~\ref{defi:Campanajets} for a precise definition.) Then a natural question is whether analogues of the existence of a section \cite[Theorem 1.1]{GHS03} and weak approximation at places of good reduction \cite[Theorem 3]{HT06} hold in the setting of Campana sections. Assuming a slightly stronger version of Campana uniruledness, we answer these questions affirmatively. 

\begin{theo}
\label{theo:mainI_intro}
Assume that $\mathbf k$ is an algebraically closed field of characteristic $0$.
Let $\pi: (\mathcal{X},\Delta_\epsilon) \to B$ be a klt Campana fibration over $\mathbf k$ such that a general fiber of $\pi$ is rationally connected and is strongly Campana uniruled (in the sense of Definition~\ref{defi:stronglycampanauniruled}).
Let $S$ be a finite set of closed points on $B$ and assume that $\underline{\pi} : \underline{\mathcal X} \to B$ satisfies weak approximation outside of $S$.
Fix a finite number of Campana jets supported on distinct fibers which are at places outside of $S$. Then this finite set of Campana jets is induced by a Campana section.
\end{theo}

In particular by \cite[Theorem 3]{HT06} this theorem applies by letting $S$ to be the set of places of bad reduction of $\underline{\pi} : \underline{\mathcal X} \to B$.

As a corollary, we obtain
\begin{coro}
\label{coro:CRC_intro}
Assume that $\mathbf k$ is an algebraically closed field of characteristic $0$.
Let $(X,\Delta_\epsilon)$ be a klt Campana orbifold over $\mathbf k$ such that $\underline{X}$ is rationally connected and $(X, \Delta_\epsilon)$ is strongly Campana uniruled.
Then $(X, \Delta_\epsilon)$ is Campana rationally connected, i.e., there is a family of Campana curves passing through two general points on $X$.
\end{coro}

Finally, we add to the list of examples where Conjecture~\ref{conj:main_intro} is known by verifying the conjecture for toric varieties.  

\begin{theo}
\label{theo:toric_intro}
Assume that $\mathbf k$ is an algebraically closed field of characteristic $0$.
Let $\underline{X}$ be a smooth projective toric variety over $\mathbf k$ and $\Delta$ be the torus-invariant boundary on $\underline{X}$.
Let $(X,\Delta_\epsilon)$ be a klt Campana orbifold. Then $(X,\Delta_\epsilon)$ is strongly Campana uniruled as well as Campana rationally connected.  
\end{theo}

\begin{rema}
In characteristic $0$, Campana noticed that for toric Campana orbifolds, Campana uniruledness easily follows by looking at toric rational curves. Here we establish something stronger so that Theorem~\ref{theo:mainI_intro} applies. We should also note that this theorem is established in positive characteristic too. See Theorem~\ref{thm:toric-uniruled-rc} for more details.
\end{rema}

As a corollary of this theorem, we deduce a version of Theorem \ref{theo:mainI_intro} when the generic fiber is toric which places no restriction on the fibers containing the jets:

\begin{coro}
\label{coro:toric}
Assume that $\mathbf k$ is an algebraically closed field of characteristic $0$.
Let $\pi: (\mathcal{X},\Delta_\epsilon) \to B$ be a klt Campana fibration over $\mathbf k$ such that 
the generic fiber $(\mathcal X_\eta, \Delta_\eta)$ is a smooth projective toric variety with the toric boundary.
Fix a finite number of Campana jets supported on distinct fibers of $\underline{\pi} : \underline{\mathcal X} \to B$. Then these finite Campana jets are induced by a Campana section.
\end{coro}

\subsection{Related works}

\subsubsection*{Weak approximation by sections}

Existence/weak approximation of sections for rationally connected fibrations has been extensively studied. As mentioned before, there is a celebrated work \cite{GHS03} showing the existence of sections in characteristic $0$. This is extended to separably rationally connected fibrations in arbitrary characteristic in \cite{deJongStarr}. Weak approximation at places of good reduction has been established in \cite{HT06} in characteristic $0$, and there are many more results in this direction, e.g., \cite{Xu12, Xu12b, Tian15, TZ18, TZ19, StarrXu, STZ22}.

\subsubsection*{$\mathbb A^1$-connectedness}

As a notion corresponding to integral points, $\mathbb A^1$-curves and $\mathbb A^1$-connectedness have been studied by various authors.
This has been first pioneered by Miyanishi and his collaborators, e.g., \cite{Miyanishi78, GM92, GMMR},
and there is an important work by Keel-McKernan \cite{KM99} on the existence of free rational curves on quasi-projective surfaces.
More recently $\mathbb A^1$-curves have been studied by the first author and Yi Zhu using log geometry and deformation theory for stable log maps (\cite{CZ15, CZ17, CZ18}) 
Surprisingly log Fano varieties with a reduced boundary do not need to satisfy $\mathbb A^1$-connectedness, e.g., the projective plane with a union of two lines as a boundary.
Moreover in characteristic $p > 0$, dlt log Fano varieties consisting of a projective space of dimension $n \geq p$ with any smooth degree $p$ boundary are at best separably $\mathbb A^1$-uniruled, and fail to be $\mathbb A^1$-separably connected in general \cite{CC21}.

\subsubsection*{Weak approximation by Campana points}

There are a few papers studying weak approximation property by Campana points in the arithmetic setting. The paper \cite{NS24} discusses a relation between weak approximation property and the Hilbert property. \cite{MNS24} is the first paper studying Brauer-Manin obstructions in the setting of semi-integral points. Finally \cite{Moerman} addresses the weak approximation property for Campana points and related notions for split toric varieties.

\subsubsection*{Campana orbifolds and orbifold rational connectedness}

Campana introduced the notion of Campana orbifolds in his studies of special manifolds, and he developed the theory of orbifold rational curves and orbifold rational connectedness in \cite{Campana07, Campana10, Campana11}. Orbifold rational connectedness is equivalent to Campana rational connectedness as defined in Definition~\ref{defi:introrc}.

\subsection{The plan of the paper}

In Section~\ref{sec:deformation}, we discuss the deformation theory of log maps and exhibit a few constructions of log gluing and log splitting which are used later. In Section~\ref{sec:deformation_logsection}, we develop the deformation theory of log sections. In Section~\ref{sec:freelogcurve}, we introduce separable uniruledness and separable connectedness by rational log curves and show that these notions are equivalent to the existence of free or very free rational log curves respectively. In Section~\ref{sec:campanamaps}, we introduce the notion of Campana maps and sections, and we also introduce the notion of Campana uniruledness and Campana rational connectedness. In Section~\ref{sec:weakapproximation}, we define Campana jets and weak approximation by Campana sections. Then we prove Theorem~\ref{theo:mainI_intro}. In Section~\ref{sec:P1}, we discuss the case of $\mathbb P^1$-fibrations and prove Conjecture~\ref{conj:main_intro} in this case. We also prove that orbifold rational connectedness is equivalent to Campana rational connectedness in Remark~\ref{rema:campanasdef}. Finally in Section~\ref{sec:toric}, we discuss the case of toric orbifolds and prove Theorem~\ref{theo:toric_intro} and Corollary~\ref{coro:toric}.

\bigskip

\noindent
{\bf Acknowledgements:}
The authors thank Brendan Hassett, Jason Starr, Yuri Tschinkel, and Chenyang Xu for answering our questions. The authors also thank Boaz Moerman for his comments on the first draft of the paper and Enhao Feng for pointing out some mistakes. Finally we would like to thank the referee for detailed suggestions which significantly improved the exposition of the paper.

Qile Chen was partially supported by NSF grant DMS-2001089.
Brian Lehmann was supported by Simons Foundation grant Award Number 851129.
Sho Tanimoto was partially supported by JST FOREST program Grant number JPMJFR212Z, by JSPS KAKENHI Grand-in-Aid (B) 23K25764, by JSPS Bilateral Joint Research Projects Grant number JPJSBP120219935, and by JSPS KAKENHI Early-Career Scientists Grant number 19K14512.

\section{Deformation theory for log maps}
\label{sec:deformation}

We work over an algebraically closed field $\mathbf{k}$.  In this section we have two goals.  First, we introduce some terminology for log schemes and stable log maps.  Second, we discuss the deformation theory of stable log maps (possibly with added constraints).

\subsection{Log maps and their stacks}\label{ss:log-map-stacks}
Let $\underline{X}$ be a smooth variety and $\Delta \subset \underline{X}$ be a strict normal crossings divisor. Denote by $\Delta = \cup_i \Delta_i$ the decomposition into smooth irreducible components.  

Let $X = (\underline{X}, \mathcal{M}_X)$ be the log scheme associated to the pair $(\underline{X},\Delta)$, where $\mathcal{M}_X$ is the sheaf of monoids over $\underline{X}$ defined by
\[
\mathcal{M}_{X}(U) := \{ f \in \cO_{\underline{X}}(U) \ | \ f|_{U \setminus \Delta} \in \cO^{\times}_{\underline{X}}(U\setminus \Delta)\}
\]
for any open subscheme $U \subset \underline{X}$. 
The log tangent bundle $T_X$ is the subsheaf of $T_{\underline{X}}$ consisting of vector fields tangent to $\Delta$. Its dual $\Omega_{X} = T_X^{\vee}$ is the {\em log cotangent bundle} consisting of differentials with at most logarithmic poles along $\Delta$.

\begin{nota}
For any log scheme (or log stack)  $Y$, we denote by $\mathcal{M}_Y$  its log structure and by $\overline{\mathcal{M}}_{Y} := \mathcal{M}_Y/\cO^*_{Y}$ the corresponding characteristic sheaf. 

For any log morphism $f \colon X \to Y$ between two log schemes or log stacks, the morphism $f^{\flat} \colon f^*\cM_Y  \to \cM_X$ (resp. $\bar{f}^{\flat} \colon f^*\overline{\cM}_Y  \to \overline{\cM}_X$) denotes the corresponding morphism on the level of log structures (resp.~characteristic sheaves). 

Any scheme $\underline{X}$ can be viewed as a log scheme with the trivial log structure $\cM_{\underline{X}} = \cO^*_{\underline{X}}$. We assume all log structures are fine and saturated, or fs for short, unless otherwise stated. 
For the basics of logarithmic structures, we refer to the foundational paper of Kato \cite{KKato89} and the comprehensive book \cite{Ogus18}. 
\end{nota}

\begin{defi}
A {\em log curve} over a log scheme $S$ consists of a pair 
\[(\pi: C \to S, \{p_1, \cdots, p_n\})\] 
such that
\begin{enumerate}
 \item The underlying pair $(\underline{\pi}, \{p_1, \cdots, p_n\})$ is a family of pre-stable curves over $\underline{S}$ with $n$ markings.
 \item $\pi$ is a proper, log smooth, and integral morphism of log schemes.
 \item If $U \subset \underline{C}$ is the smooth locus of $\underline{\pi}$ then $\overline{\mathcal{M}}_{C}|_{U} \cong \underline{\pi}^*\overline{\mathcal{M}}_{S}\oplus \bigoplus^{n}_{k=1} p_{k*}\mathbb{N}_{\underline{S}}$.   
\end{enumerate}
Here $\mathbb{N}_{\underline{S}}$ denotes the constant sheaf on $\underline{S}$ with coefficients $\mathbb{N}$. 
\end{defi}

Let $S$ be a log scheme. A {\em log map} over $S$ 
is a morphism of log schemes $f \colon C \to X$ such that $C \to S$ is a family of log curves. In particular, the underlying family $\underline{C} \to \underline{S}$ obtained by removing all log structures is  a family of pre-stable curves. It is called {\em stable} if the corresponding underlying morphism $\underline{f}$ is stable in the usual sense. 
A log map $f$ is said to be {\em non-degenerate}
if $S$ is a log point with the trivial log structure. In this case, $\underline{C}$ is a smooth irreducible curve and $\mathcal{M}_C$ is the divisorial log structure coming from the markings.  Furthermore in this case we can conclude that $f(C) \not\subset \Delta$.

The theory of log maps of Abramovich--Chen--Gross--Siebert \cite{Ch14, AbCh14, GrSi13} is one of the main tools in this paper. 
Let $f \colon C \to X$ be a non-degenerate log map and let $p_k \in C$ be the $k$-th marked point. Let $c_{k,i}$ be the order of tangency of $f$ with respect to $\Delta_i$ at $p_k$.  The collection of non-negative integers $\mathbf{c}_k = (c_{k,i})_i$ is called the {\em contact order} at $p_k$. Log geometry allows us to further extend this definition to all log maps, not only the non-degenerate ones.

We recall the definition of contact orders at a node as described by \cite[\S 3.2]{Ch14} and \cite[\S 4.1]{AbCh14}. Let $f \colon C \to X$ be a log map over a log point $S$, and $p \in C$ be a node with image $f(p) = x$. Define the relative characteristic sheaf
$
\ocM_{C/S} := \ocM_{C}/\ocM_{S}.
$
Recall that $\ocM_{C/S}|_{p} \cong \mathbb{Z}$ where this isomorphism depends on a choice of sign. Consider the composition 
\[
\xymatrix{
u_p \colon \ocM_{X}|_{x} \ar[r]^{\bar{f}^{\flat}|_p} & \ocM_{C}|_{p} \ar[r] & \ocM_{C/S}|_{p} \cong \mathbb{Z}. 
}
\] 
Suppose $J  = \{i \, | \, x \in \Delta_i\}$. 
Then $\ocM_{X}|_{x} \cong \mathbb{N}^{|J|}$ with the generator $\delta_i \in \ocM_{X}|_{x}$ corresponding to $\Delta_i$. Write $c_{p,i} := u_p(\delta_i)$ if $i \in J$ and $c_{p,i} = 0$ otherwise. As a morphism of monoids, $u_p$ is uniquely determined by the collection of integers $(c_{p,i})_i$.  The contact order at the node $p$ is the collection $(c_{p,i})_i$.
Note that the contact order at a node depends on the choice of sign in the isomorphism $\ocM_{C/S}|_{p} \cong \mathbb{Z}$. However, this difference will not be important in this paper; we will only care about the divisibility properties of the contact order.

\begin{nota}\label{nota:contact-order}
Consider a contact order $\mathbf{c}_k = (c_{k,i})_i$ of a marking or a node. We say that a contact order is {\em positive} if all the entries are non-negative and at least one entry is not zero.  A marking or a node is a {\em contact marking} if its contact order is positive; otherwise, it is called a {\em non-contact marking}. 
For an integer $m \in \mathbb{N}$, we wrote $\mathbf{c}_k \geq m$ if $c_{k,i} \geq m $ for all $i$. We write $\ch \mathbf{k} \mid \mathbf{c}_k$ if $\ch \mathbf{k}  \mid c_{k,i}$ for all $i$, and write $\ch \mathbf{k} \nmid \mathbf{c}_k$ otherwise. 

For later use, denote by $\underline{X}_{\mathbf{c}_{k}} := \cap_{c_{k,i} \neq 0} \Delta_i$ with the trivial log structure if $\mathbf{c}_k \neq 0$, and $\underline{X}_{\mathbf{c}_k} = \underline{X}$ if $\mathbf{c}_k = 0$.  Further denote by $\underline{X}^{\circ}_{\mathbf{c}_{k}} \subset \underline{X}_{\mathbf{c}_{k}}$ the maximal open dense locus
such that $\underline{X}^{\circ}_{\mathbf{c}_{k}} \cap \Delta_i = \emptyset $ for any index $i$ such that $c_{k,i} = 0$. 
\end{nota}

The discrete data of a stable log map to $X$ is the triple
\begin{equation}\label{eq:discrete-data}
(g, \varsigma = \{\mathbf{c}_k\}_{k=1}^{|\varsigma|}, \beta)
\end{equation}
where $g$ denotes the genus of the domain curve, $\beta$ is a curve class on $\underline{X}$, $|\varsigma|$ is the number of markings, and $\mathbf{c}_k$ is the {\em contact order} at the $k$-th marked point.

Let $\mathscr{M}_{g, \varsigma}(X,\beta)$ be the category of stable log maps to $X$ with the discrete data \eqref{eq:discrete-data} fibered over the category of log schemes. It was proved in  \cite{Ch14, AbCh14, GrSi13, Wise16} that $\mathscr{M}_{g, \varsigma}(X,\beta)$ is represented by a log algebraic stack, i.e.~it is an algebraic stack equipped with a fine and saturated log structure. Further assuming characteristic zero, then $\mathscr{M}_{g, \varsigma}(X,\beta)$ is a proper, log Deligne-Mumford stack. 
Let $\mathscr{M}^{\circ}_{g, \varsigma}(X,\beta) \subset \mathscr{M}_{g, \varsigma}(X,\beta)$ be the open sub-stack with the trivial log structure. Then  $\mathscr{M}^{\circ}_{g, \varsigma}(X,\beta)$ is the stack parametrizing non-degenerate stable log maps with the discrete data assigned in \eqref{eq:discrete-data}.

\subsection{Deformations of log maps}\label{ss:deformation}

For any log stack $\mathfrak{M}$, denote by $\mathbf{Log}_{\mathfrak{M}}$ Olsson's log stack constructed in \cite{Olsson03}: to any underlying morphism $\underline{S} \to \underline{\mathfrak{M}}$ the stack associates the category of morphisms of log stacks $S \to \mathfrak{M}$. The universal log structure of $\mathbf{Log}_{\mathfrak{M}}$ defines the tautological morphism of log stacks $\mathbf{Log}_{\mathfrak{M}} \to \mathfrak{M}$. Consider 
\[
\mathfrak{M}^{\log}_{g, |\varsigma|} := \mathbf{Log}_{\mathfrak{M}_{g, |\varsigma|}}.
\]
where $\mathfrak{M}_{g, |\varsigma|}$ is the moduli stack of pre-stable curves equipped with the canonical log structure \cite{FKato00, Olsson07}.
The stack $\mathfrak{M}^{\log}_{g, |\varsigma|}$ is log smooth. Hence the locus in $\mathfrak{M}^{\log}_{g, |\varsigma|}$ with the trivial log structure is open dense and maps isomorphically via the tautological morphism to the open substack of $\mathfrak{M}_{g, |\varsigma|}$ parametrizing smooth curves. 
In particular, $\mathfrak{M}^{\log}_{g, |\varsigma|}$ is reduced and irreducible of dimension  
\[
\dim \mathfrak{M}^{\log}_{g, |\varsigma|} = 3g-3 + |\varsigma|. 
\]

Let $\mathfrak{C}_{g, |\varsigma|} \to \mathfrak{M}_{g, |\varsigma|}$ be the universal log curve. Consider the pull-back of log stacks
\[
\mathfrak{C}^{\log}_{g, |\varsigma|} := \mathfrak{C}_{g, |\varsigma|} \times_{\mathfrak{M}_{g, |\varsigma|}}\mathfrak{M}^{\log}_{g, |\varsigma|} \to \mathfrak{M}^{\log}_{g, |\varsigma|}.
\]
It is a universal family in the following sense.
Let $C \to S$ be a log curve of genus $g$ with $|\varsigma|$ markings.
 Then there is a natural (not necessarily strict) morphism $S \to \mathfrak{M}_{g, |\varsigma|}$ such that $C \to S$ is the pull-back $C = \mathfrak{C}_{g, |\varsigma|} \times_{\mathfrak{M}_{g, |\varsigma|}}S \to S$. The functoriality of Olsson's log stack implies that $S \to \mathfrak{M}_{g, |\varsigma|}$ factors through a unique strict morphism $S \to \mathfrak{M}^{\log}_{g, |\varsigma|}$ such that the family $C \to S$ is the pull-back $C = \mathfrak{C}^{\log}_{g, |\varsigma|} \times_{\mathfrak{M}^{\log}_{g, |\varsigma|}}S \to S$. 
 In particular, we obtain a tautological morphism 
\begin{equation}\label{eq:forget-log-map}
\mathscr{M}_{g, \varsigma}(X,\beta) \longrightarrow \mathfrak{M}^{\log}_{g, |\varsigma|}, \ \ [f\colon C/S \to X] \mapsto [C/S]
\end{equation}
which assigns to a log map its domain log curve.

\subsubsection{Deformations of log maps relative to $\mathfrak{M}^{\log}_{g,|\varsigma|}$}
Suppose $S$ is a geometric log point of $\mathscr{M}_{g, \varsigma}(X,\beta)$ corresponding to a log map $f: C \to X$. The first-order deformations and obstructions of the log map
relative to $\mathfrak{M}^{\log}_{g,|\varsigma|}$ are controlled by $H^0(f^*T_X)$ and $H^1(f^*T_X)$ respectively (\cite[\S 4]{ACGS21}). In particular, the expected relative dimension at $[f]$ is 
\[
\exp \dim_{[f]} \left( \mathscr{M}_{g, \varsigma}(X,\beta) \Big/ \mathfrak{M}^{\log}_{g, |\varsigma|} \right) = \chi(f^*T_X),
\]
yielding the expected dimension
\begin{equation}\label{eq:expected-dim-map}
\exp \dim_{[f]} \left( \mathscr{M}_{g, \varsigma}(X,\beta) \right) = \chi(f^*T_X) + 3g-3 + |\varsigma|. 
\end{equation}

\begin{exam}[Toric example]
Let $X$ be the log scheme associated to a toric variety with its toric boundary. Then $T_X \cong \cO^{\oplus \dim X}$ by \cite[Theorem 8.1.1]{CLS11}.  For any genus zero log map $f \colon \cC \to X$ over a geometric log point $S$, we have $H^1(f^*T_X) = 0$. This implies that  the tautological morphism \eqref{eq:forget-log-map} is log smooth.  Consequently $\mathscr{M}_{0, \varsigma}(X,\beta)$ is also log smooth of dimension  
\[
\dim \left( \mathscr{M}_{0, \varsigma}(X,\beta) \right) = \dim X + |\varsigma| - 3. 
\]
\end{exam}

\subsubsection{Deformations of log maps relative to $\mathbf{Log}$}
The morphism of log stacks $\mathfrak{M}_{g,|\varsigma|} \to \Spec \, \mathbf{k}$ induces a tautological morphism 
\[
\mathfrak{M}^{\log}_{g, |\varsigma|} \longrightarrow \mathbf{Log} := \mathbf{Log}_{\Spec \mathbf{k}}, \qquad [C\to S] \mapsto [S], 
\]
hence a tautological morphism by composing \eqref{eq:forget-log-map}
\begin{equation}\label{eq:to-log}
\mathscr{M}_{g, \varsigma}(X,\beta) \longrightarrow \mathbf{Log}, \ \ [f\colon C/S \to X] \mapsto [S]
\end{equation}

For a log map $f \colon C \to X$ over $S$, consider the complex $N_{f}$ defined by the distinguished triangle
\begin{equation}\label{eq:deformation-complex}
T_{C/S} \stackrel{d f}{\longrightarrow} T_{X} \longrightarrow N_{f}  \stackrel{[1]}{\longrightarrow}
\end{equation}
Suppose $S$ is a geometric log point. The first-order deformations and obstructions of $[f]$
relative to $\mathbf{Log}$ are controlled by $H^0(N_{f})$ and $H^1(N_{f})$ respectively. In general $N_{f}$ is a complex rather than a sheaf. We call $N_f$ the {\em normal complex} of $f$. %

\subsection{The deformation theory of log immersions}\label{ss:immersion-deformation}

We describe a situation where the complex $N_f$ is represented by a vector bundle. This will allow us to effectively compute deformations of log maps in many examples.  

\begin{defi}\label{defi:log-immersion}
Let $f \colon C \to X$ be a log map over a geometric log point $S$.  We say that $f$ is a {\em log immersion} if 
\begin{enumerate}
 \item For any node or marking $p \in C$ with contact order $\mathbf{c}_p$, we have $\mathbf{c}_p \neq 0$ and $\ch\mathbf{k} \nmid \mathbf{c}_p$. 
 \item The underlying morphism $\underline{f} \colon \underline{C} \to \underline{X}$ is an immersion away from nodes and markings.  
\end{enumerate}
\end{defi}

Condition $(1)$ implies the images of all nodes and markings are necessarily in $\Delta$. 
The following is a generalization of \cite[Lemma 4.12]{CZ17} with a similar proof. We provide a detailed proof for completeness. 

\begin{lemm}\label{lem:normal-bundle}
Suppose $f \colon C \to X$ is a log immersion. Then $d f \colon T_{C/S} \to f^*T_X$ is injective as a subbundle. Hence the cokernel  $N_f := \cok{d f}$ is a vector bundle on $C$. 
\end{lemm}
\begin{proof}
It suffices to prove dually that $f^* \colon f^*\Omega_X \to \Omega_{C/S}$ is surjective and the kernel is locally free.
We verify this locally around a closed point $x \in C$. 

\smallskip
\noindent
{\bf Case 1: Smooth non-marked points.} 
Suppose $x$ is a smooth non-marked point. Let $Z \subset C$ be the irreducible component containing $x$. For a subset $J \subset \{i\}$, denote by $\underline{X}_{J}^\circ := (\cap_{j\in J}\Delta_j) \setminus \cup_{i\not\in J}\Delta_i$, and set $\underline{X}_{\emptyset}^\circ = \underline{X} \setminus \Delta$ for $J = \emptyset$.  
Thus $\{\underline{X}_{J}^\circ \ | \ \underline{X}_{J}^\circ \neq \emptyset \}$ is a stratification of $\underline{X}$ such that $\ocM_{X}|_{\underline{X}_J^\circ}$ is the trivial monoid $\mathbb{N}^{|J|}$. Let $Z^{\circ} \subset Z$ be the open subscheme obtained by removing all nodes and markings. It contains $x$.  Then there is a unique stratum $\underline{X}_{J}^\circ \neq \emptyset$ such that the underlying morphism $\underline{f}|_{\underline{Z}^{\circ}} \colon \underline{Z}^{\circ} \to \underline{X}$ factors through $\underline{X}_{J}^\circ$.
Thus the induced $\underline{Z}^{\circ} \to \underline{X}_{J}^\circ$ is an immersion in the usual sense. We have a commutative diagram
\[
\xymatrix{
0 \ar[r] &  \Omega_{\underline{X}_J^\circ}|_{\underline{Z}^{\circ}} \ar[d] \ar[r] & \Omega_{X}|_{\underline{Z}^{\circ}} \ar[d]^{f^*|_{\underline{Z}^{\circ}}} \ar[r] & \cO_{\underline{Z}^{\circ}}^{\oplus |J|} \ar[r] & 0 \\
& \Omega_{\underline{Z}^{\circ}} \ar[r]^{\cong} & \Omega_{C/S}|_{\underline{Z}^{\circ}} & & 
}
\]
where the top is an exact sequence by \cite[Lemma 4.13]{CZ17}
and $\Omega_{\underline{X}_{J}}$ is the cotangent bundle of the underlying scheme. 
The immersion $\underline{Z}^{\circ} \to \underline{X}_{J}^\circ$ implies the surjectivity of $f^*|_{\underline{Z}^{\circ}}$  as needed. Moreover the kernel is locally free in a neighborhood of $x$ because $\Omega_{X}|_{\underline{Z}^{\circ}}$ is torsion free and $x$ is a smooth point of $\underline{C}$.

\smallskip
\noindent
{\bf Case 2: Marked points.}
Let $x$ be a marked point, defined by a local coordinate $z$. Thus the fiber $\omega_{C/S}|_{x}$ is generated by a section $\frac{d z}{z}$. 
Let $\underline{X}_{J}$ 
be the unique stratum containing the image $f(x)$. Thus the fiber $\Omega_{X}|_{x}$ contains a set of linearly independent vectors $\{\frac{d \delta_j}{\delta_j}|_{x} \ | \ j \in J \}$, where $\delta_j$ are defining equations of $\Delta_j$ locally around $f(x)$. Let $\mathbf{c}_x$ be the contact order at $x$; we write it in the form $\mathbf{c}_{x} = (c_{x,j})_{j \in J}$ 
where the contact order at $x$ along $\Delta_j$ is given by $c_{x,j} \in \mathbb{N}$. 
By assumption, there exists $j \in J$ such that $c_{x,j} > 0$ and $\ch \mathbf{k} \nmid c_{x,j}$. Thus, by an appropriate choice of coordinates, we have $f^*(\frac{d \delta_j}{\delta_j}|_{x}) = c_{x,j}\cdot\frac{d z}{z} \neq 0$. Hence $f^*|_{x}$ is surjective. As $f^*$ is a morphism of coherent sheaves, this implies the surjectivity of $f^*$ in a neighborhood of $x$.  
Again the kernel is locally free because $x$ is a smooth point of $\underline{C}$.

\smallskip
\noindent
{\bf Case 3: Nodal points.} Suppose $x$ is a node. Similar to the case of marked points, it suffices to prove the surjectivity of $f^*|_{x}$ and local freeness of the kernel on a neighborhood of $x$.  Suppose $x$ is a node joining two branches $Z_1$ and $Z_2$ with local coordinates $z_1$ and $z_2$. Then $\omega_{C/S}|_{x}$ is generated by a local section $\frac{d z_1}{z_1} = -\frac{d z_2}{z_2}$.
Again let $\underline{X}_{J}$ be the unique stratum containing the image $f(x)$. Thus the fiber $\Omega_{X}|_{x}$ contains a set of linearly independent vectors $\{\frac{d \delta_j}{\delta_j}|_{x} \ | \ j \in J \}$ where each $\delta_j$ is the defining equation of $\Delta_j$ locally around $f(x)$. 
Let $\mathbf{c}_x$ be the contact order at the node $x$. (The definition was given just before Notation~\ref{nota:contact-order}.)
This means that for each $j \in J$, one of the two branches $Z_i$ has contact order 
$c_{x,j} \in \mathbb{N}$ along $\Delta_j$ at $x$. 
The assumption on contact orders implies that there is a $j \in J$ and an $i \in \{1,2\}$ such that $Z_{i}$ has contact order $c_{x,j} > 0$ along $\Delta_{j}$ such that $\ch \mathbf{k} \nmid c_{x,j}$. Similarly, by an appropriate choice of coordinates, we have $f^*(\frac{d \delta_j}{\delta_j}|_{x}) = c_{x,j}\cdot\frac{d z_i}{z_i} \neq 0$. This implies the desired surjectivity of $f^*|_{x}$. 
Furthermore, an element $\sum_{j \in J} a_j \frac{d \delta_j}{\delta_j}|_{x}$ lies in the kernel of $f^*|_x$ if it satisfies the linear relation $\sum_{j \in J}c_{x,j}a_{j} = 0$. Thus the kernel has codimension $1$ so that it is locally free at $x$.
\end{proof}

Consider a log map $f \colon C \to X$ over a log point $S$. There is a natural {\em associated log map} $\bar{f} \colon \bar{C} \to X$ over $S$ such that the composition $C \stackrel{}{\lra} \bar{C} \stackrel{}{\lra} X$ is $f$ where the first arrow is the forgetful morphism removing all markings with the zero contact order. In particular, $\bar{f}$ is a log map with only contact markings. 

\begin{coro}\label{cor:torsion-normal-complex}
Notations as above, further assume that $\bar{f}$ is a log immersion and let $P$ denote the set of non-contact markings of $f$. Then $N_{f}$ is a sheaf with torsion-free part $N_{f}^{tf} = N_{\bar{f}}$ and torsion $N_{f}^{tor} = \oplus_{k \in P} N_{p_k/\underline{C}}$. 
\end{coro}
\begin{proof}
Consider the following diagram whose rows are distinguished triangles:
\[
\xymatrix{
T_{C/S} \ar[r] \ar[d] &  f^*T_{X} \ar[r] \ar[d]_{=} & N_{f} \ar[r]^{[1]} \ar[d] &  \\
T_{\bar{C}/S} \ar[r] & \bar{f}^*T_{X} \ar[r] & N_{\bar{f}} \ar[r]^{[1]} & 
}
\]
The injectivity of $T_{C/S} \to T_{\bar{C}/S}$ and Lemma \ref{lem:normal-bundle} imply that $T_{C/S} \to  f^*T_{X}$ is injective. Hence $N_{f}$ is a sheaf. 
By Lemma \ref{lem:normal-bundle} again, $N_{\bar{f}}$ is locally free.  The kernel-cokernel sequence shows that the kernel of $N_{f} \to N_{\bar{f}}$ can be identified with the cokernel of the leftmost map, namely $\oplus_{k \in P} N_{p_k/\underline{C}}$. This implies our assertion.
\end{proof}

\subsection{Deformations of log maps with point constraints}\label{ss:point-constraints-deformation}
For a subset $P \subset \{1, \cdots, |\varsigma|\}$ and a collection of points $q_k \in X^{\circ} := X \setminus \Delta$ for $k \in P$, we introduce the point constraint 
\[
f_P := \{ p_k \mapsto q_k \ | \ k \in P \}
\] 
sending the $k$-th marking to $q_k$. A log map $f$ satisfies the point constraint $f_P$ iff $f(p_k) = q_k$ for $k \in P$. As $q_k \in X^{\circ}$, we necessarily have $\mathbf{c}_k = 0$.

For each $k \in P$, consider the evaluation morphism 
\[
\ev_k \colon \mathscr{M}_{g, \varsigma}(X,\beta) \longrightarrow \underline{X}, \ \ [f] \mapsto f(p_k)
\]
induced by the $k$-th marking. The moduli of stable log maps satisfying the point constraint $f_P$ is 
\begin{equation}\label{eq:moduli-point-constraints}
\mathscr{M}_{g, \varsigma}(X,\beta, f_P) := \mathscr{M}_{g, \varsigma}(X,\beta)\times_{\prod_P \underline{X}} \prod_{k \in P} q_k
\end{equation}
with $\prod_{k \in P} \ev_k \colon \mathscr{M}_{g, \varsigma}(X,\beta) \to \prod_P \underline{X}$. 

Consider a log map $f \colon C \to X$ over $S$ satisfying the point constraint $f_P$. As $\mathbf{c}_k = 0$ for $k \in P$, by removing $\{p_k \ | \ k \in P\}$ from the set of markings, $f$ induces an associated log map $\bar{f} \colon \bar{C} \to X$ over $S$ such that $f$ is given by the composition $C \to \bar{C} \to X$. Consider the twisted normal complex 
\begin{equation}\label{eq:twisted-normal-complex}
N_{f, f_P} := N_{\bar{f}}\left(-\sum_{k\in P}p_k \right)
\end{equation}
Assuming that $S$ is a geometric log point, then the first order deformations and obstructions of $f$ relative to $\mathbf{Log}$ are given by $H^0(N_{f, f_P})$ and $H^1(N_{f, f_P})$ respectively (\cite[\S 4]{ACGS21}).

\subsection{Deformation of log maps relative to codimension $1$ boundary}

The condition Definition \ref{defi:log-immersion}.(1) leads to a very intuitive geometric picture around contact markings as observed below. 

\begin{lemm}\label{lemm:N-at-marking}
Let $f \colon C \to X$ be a log map (not necessarily a log immersion) over a geometric log point $S$. Suppose $p \in C$ is a marking satisfying $\ch\mathbf{k} \nmid \mathbf{c}_p$. Further assume that the image $x = f(p)$ is contained in a unique boundary component $\Delta_i \subset \Delta$. Then there is a natural isomorphism 
\[
T_{{\Delta}_i}|_{x} \cong N_{f}|_{p}. 
\]
\end{lemm}
\begin{proof}
We have the following commutative diagram of solid arrows at $p$: 
\[
\xymatrix{
0 \ar[d]&& 0\ar[d] && && \\
\cO_{p} \ar[d]_{\cong} \ar[rr]^{\cong} && \cO_{p} \ar[d] \ar@{-->}[rrd]^0 && &&  \\
T_{C/S}|_p \ar[d] \ar[rr]^{df} && T_{X}|_{p} \ar[rr] \ar[d] && N_{f}|_{p} \ar[rr]^{[1]} && \\
0 \ar[d] \ar[rr] && T_{{\Delta_i}}|_{p} \ar[d] \ar@{-->}[rru]_{\cong} && && \\
0  && 0 && &&
}
\]
where the middle row is the pull-back of \eqref{eq:deformation-complex} to $p$ and the two columns are obtained by pulling back corresponding exact sequences in \cite[Lemma 4.1]{CZ14}. The same calculation as in Lemma \ref{lem:normal-bundle}, case 2 shows that $df$ is non-zero, inducing the isomorphism $\cO_p \cong \cO_p$. Consequently, we obtain the two dashed arrows, finishing the proof. 
\end{proof}

Consider the moduli $\mathscr{M}_{g, \varsigma}(X,\beta, f_P)$ as in \eqref{eq:moduli-point-constraints}. Let $p_k$ be a contact marking such that $k\not \in P$. Further assume there is an $i'$ such that $\ch \mathbf{k} \nmid c_{k,i'} \neq 0$, and  $c_{k,i} = 0$ for all $i \neq i'$.
This induces an evaluation morphism 
\[
\ev_k \colon \mathscr{M}_{g, \varsigma}(X,\beta, f_P) \to \Delta_{i'}.
\]
We are interested in the local structure of this morphism along the stratum $\Delta^{\circ}_{i'} = \Delta_{i'}\setminus \cup_{i\neq i'}\Delta_i$. 

\begin{prop}\label{prop:deformation-relative-boundary}
Consider a log map $[f \colon C \to X] \in \mathscr{M}^{\circ}_{g, \varsigma}(X,\beta, f_P)(S)$ over a geometric log point $S$ such that $x = f(p_k) \in \Delta^{\circ}_{i'}$. Then there is a natural morphism 
\[
d \ev_k|_{[f]} \colon H^{0}(N_{f, f_{P}}) \to T_{\Delta_{i'}}|_{x}. 
\]
It is surjective if $H^{1}(N_{f,f_P}(- p_k)) = 0$. 
\end{prop}
\begin{proof}
Consider the exact sequence
\[
0 \to \cO_{C}(-p_k) \to \cO_{C} \to \cO_{p_k} \to 0
\]
Tensoring with $N_{f, f_P}$, we obtain a distinguished triangle
\[
N_{f,f_P}(- p_k) \to N_{f,f_P} \to N_{f, f_{P}}|_{p_k} \stackrel{[1]}{\to}
\]
Taking the long exact sequence of cohomology, we have 
\[
H^0(N_{f,f_P}(- p_k)) \to H^0(N_{f,f_P}) \to H^0(N_{f, f_{P}}|_{p_k}) \to H^1(N_{f,f_P}(- p_k)) 
\]
Then $d \ev^{\circ}_k|_{[f]} $ is given by the composition 
\[
H^0(N_{f,f_P}) \to H^0(N_{f, f_{P}}|_{p_k}) \cong T_{{\Delta}_i}|_{x} 
\]
with the isomorphism given by Lemma \ref{lemm:N-at-marking}.
When $ H^1(N_{f,f_P}(- p_k))= 0$ surjectivity follows from the long exact sequence. 
\end{proof}

\subsection{Log gluing along markings}
\label{subsec:gluing}

We provide a gluing construction of log maps that can increase contact orders at a given marking. Similar constructions for $\mathbb{A}^1$-curves in different settings were provided in \cite[\S 4]{CZ15} and \cite[\S 4]{CZ18}. 

Suppose that for each $j=1,2$, we have a log map $f_j \colon C_j \to X$  over a geometric log point $S_i$, and a marking $p_j \in C_j$ with contact order $\mathbf{c}_{p_j}$ contained in an irreducible component $Z_j \subset C_j$ such that
\[
\underline{S}_1 = \underline{S}_2, \qquad \mbox{and} \qquad f_{1}(p_1) = f_2(p_2). 
\]
Set $x = f_i(p_i) \in X$. We further assume that $f_j(Z_j) \not\subset \Delta_{\mathrm{red}}$ for $j=1,2$. For later use, let $M_j$ denote the set of markings on $C_j$. 

In the following subsections we describe how to construct a log stable map $f: C \to X$ that restricts to the maps $f_{1},f_{2}$ along certain irreducible components of $C$.  The construction depends on how the intersection point $x$ interacts with the irreducible components of the boundary divisor.  We only describe the construction when $x$ lies on at most one irreducible component.

\subsubsection{Gluing along non-degenerate point}

First, we assume that $x \in X^{\circ}$. In particular $\mathbf{c}_{p_1} = \mathbf{c}_{p_2} = 0$. To glue $f_1$ and $f_2$ along $p_1$ and $p_2$, we first construct the glued underlying pre-stable map $\underline{f} \colon \underline{C} \to  \underline{X}$ such that $\underline{C} = \underline{C}_1 \cup_{p_1= p_2} \underline{C_2}$ is obtained by identifying $p_1$ and $p_2$, and $\underline{f}|_{\underline{C}_j} = \underline{f}_j$ for each $j$.

Let $M$ be the set of markings on the new prestable curve $\underline{C} \to \underline{S}$. It is given by
\[
M := (M_1\setminus \{p_1\}) \sqcup (M_2\setminus \{p_2\}).
\] 
Note that $p_1, p_2$ are glued to a node, denoted by $p$. Let $C^{\sharp} \to S^{\sharp}$ be the canonical log curve over $\underline{C} \to \underline{S}$ taking into account all markings. Similarly, let $C_j^{\sharp} \to S_j^{\sharp}$ be the canonical log curve over $\underline{C}_j \to \underline{S}_j$. Over the same underlying point $\underline{S}$, we have a splitting
\[
\cM_{S^{\sharp}} = \cM_{S_1^{\sharp}} \oplus_{\cO^{\times}} \cN_{p} \oplus_{\cO^{\times}} \cM_{S_2^{\sharp}}
\]
where $\cN_{p}$ is the canonical log structure smoothing the node $p$.

Recall that the stable log map $C_j \to S_j$ is equivalent to a morphism of log structures $\cM_{S_j^{\sharp}} \to \cM_{S_j} $, i.e., this means that $C_j \to S_j$ is isomorphic to $C_j^\sharp \times_{S_j^\sharp} S_j \to S_j$.
We then define the log point $S$ by setting 
\[
\cM_S := \cM_{S_1} \oplus_{\cO^{\times}} \cN_{p} \oplus_{\cO^{\times}} \cM_{S_2}.
\]
The morphism of log structures $\cM_{S^{\sharp}} \to \cM_S$ induced by $\cM_{S_j^{\sharp}} \to \cM_{S_j} $ and the identity $\cN_p \to \cN_p$ defines a morphism of log points $S \to S^{\sharp}$. We obtain the log curve $C := C^{\sharp}\times_{S^{\sharp}}S \to S$. 

Finally, $\underline{f}$ lifts to a stable log map $f \colon C \to X$ such that the restriction $f|_{\underline{C}_j}$ is induced by $f_j$ naturally. 
Note that for every marking $p' \in C$ its contact order is the same as the contact order of the corresponding marking from $f_j$.

\subsubsection{Gluing along markings in codimension $1$ boundary strata}\label{sss:gluing-along-stratum}

Now we assume that there is an index $i$ such that $x \in \Delta_{i'}$ iff $i' = i$. In particular, $x$ is contained in the unique codimension $1$ boundary stratum of $\Delta_i$.  
For simplicity we will write $c_{p_{1}}$ and $c_{p_{2}}$ for the unique non-zero element in the contact orders $\mathbf{c}_{p_{1}}$ and $\mathbf{c}_{p_{2}}$.

\smallskip
\noindent
{\bf Step 1. The underlying pre-stable map.}
To glue $f_1$ and $f_2$, we introduce a smooth rational curve $\underline{L}$ with three distinct special points $p'_1, p'_2, p$. We form the underlying pre-stable curve over $\underline{S}$
\[
\underline{C} = \underline{C}_1 \cup_{p_1 = p'_1} \underline{L} \cup_{p'_2 = p_2} \underline{C}_2
\]
by the corresponding identifications. Denote by $q_j$ the node obtained by gluing $p_j$ and $p'_j$. 
The set of markings on $\underline{C}$ is given by 
\[
M := (M_1\setminus \{p_1\}) \sqcup \{p\} \sqcup (M_2\setminus \{p_2\}).
\]
The glued underlying pre-stable map $\underline{f} \colon \underline{C} \to \underline{X}$ is defined by 
\[
\underline{f}|_{\underline{C}_j} := \underline{f}_j \qquad \mbox{and} \qquad \underline{f}(\underline{L}) = x.
\]

\smallskip
\noindent
{\bf Step 2. The domain log curve.}
Next, we construct a stable log map $C \to S$ over the underlying curve $\underline{C} \to \underline{S}$ as follows. 
Let $C^{\sharp} \to S^{\sharp}$ be the canonical log curve over $\underline{C} \to \underline{S}$ with the set of markings $M$. There is a splitting
\[
\cM_{S^{\sharp}} = \cM_{S_1^{\sharp}} \oplus_{\cO^{\times}} \cN_{q_1}   \oplus_{\cO^{\times}} \cN_{q_2} \oplus_{\cO^{\times}} \cM_{S_2^{\sharp}}
\]
where $\cN_{q_j}$ is the canonical log structure on $\underline{S}$ smoothing the node $q_j$, and $\cM_{S_j^{\sharp}}$ is the canonical log structure associated to the underlying curve $\underline{C}_j \to \underline{S}$ as above. 

Since $\underline{S}$ is a geometric point, we have a (non-canonical) splitting 
\[
\cN_{q_j} \cong \mathbb{N}\times \cO^{\times}_{\underline{S}}.
\]
Define a new log structure $\cN := \mathbb{N}\times \cO^{\times}_{\underline{S}}$ with the structure arrow $\alpha \colon \cN \to \cO_{\underline{S}}$ such that $\alpha(n,u) = 0$ if $n > 0$ and $\alpha(n,u) = u$ if $n=0$.  Set $\ell = \operatorname{lcm}(c_{p_1}, c_{p_2})$\footnote{In case $c_{p_1} = c_{p_2} = 0$, we may set $\ell = 1$ and the following construction still works. }. We fix a morphism of log structures for each $j$
\[
\cN_{q_j} \longrightarrow \cN,\qquad (n,u) \mapsto  (n\cdot \ell/c_{p_j}, u). 
\]
We then define $S$ to be the geometric log point with log structure 
\[
\cM_{S} := \cM_{S_1} \oplus_{\cO^{\times}} \cN \oplus_{\cO^{\times}} \cM_{S_2}.
\]
This leads to a morphism of log structures
\[
\cM_{S^{\sharp}} \to \cM_{S}
\]
induced by $\cM_{S_j^{\sharp}} \to \cM_{S_j} $ and the fixed morphism $\cN_{q_j} \to \cN$ above, hence a morphism of log points $S \to S^{\sharp}$. This defines a morphism of log points $S \to S^{\sharp}$, hence a log curve 
\[
C := C^{\sharp}\times_{S^{\sharp}}S \to S.
\]

\smallskip
\noindent
{\bf Step 3. The log map defined over $C \setminus L^\circ$.}  
Denote by $L \subset C$ the closed strict subscheme over $\underline{L}$ and let $L^\circ = L\setminus \{q_1, q_2\}$.  Then we have $C \setminus L = C^{\circ}_1 \sqcup C^{\circ}_2$ where $\underline{C}^{\circ}_j = \underline{C}_j \setminus \{p_j\}$. The construction of $C \to S$ implies that 
\[
\cM_{C^{\circ}_j} = \cM_{C_j}|_{\underline{C}^{\circ}_j}\oplus_{\cO^{\times}} \cN \oplus_{\cO^{\times}} \cM_{S_{j'}},
\]
where $j' \neq j$. Thus we naturally define
\[
f|_{C^{\circ}_j} := f_j|_{\underline{C}^{\circ}_j} \colon C^{\circ}_j \to X
\]
for $j=1,2$. 

Further recall that $f_{j}(Z_j) \not\subset \Delta_{red}$, i.e. $f|_{Z_j}$ intersects the boundary $\Delta$ properly at $q_j$. We observe that $f|_{C^{\circ}_j}$ extends uniquely across the node $q_j$ to yield a morphism 
\begin{equation}\label{eq:log-map-away-from-linker}
f|_{\overline{C^{\circ}_1 \sqcup C^{\circ}_2}} \colon \overline{C^{\circ}_1 \sqcup C^{\circ}_2} \to X
\end{equation}
where $\overline{C^{\circ}_1 \sqcup C^{\circ}_2} = C \setminus \left( L^\circ \right)$.  
To see this, let $s_j, t_j$ be the local coordinates around $q_j$ on the components $Z_j$ and $L$ respectively.   They each lift uniquely to a local section of $\cM_C$. 

On the other hand, let $\delta$ be a local section of $\cO_{X}$ around $x$ which is a local defining equation of $\Delta_i$. It lifts to a unique local section in $\cM_X$ around $x$, denoted by $\delta$ again. Consequently, the pull-back log structure $\underline{f}^*\cM_{X}$ around the point $q_j$ is generated by a unique section, denoted again by $\delta$. On the level of underlying morphisms, we have 
\[
\underline{f}^*(\delta) = s_j^{c_{p_j}}
\]
by choosing coordinates properly. As morphisms of log schemes are required to be compatible with their underlying structures, on the level of log structures locally around $q_j$ we have 
\begin{equation}\label{eq:gluing-node-map}
f^{\flat}(\delta) = c_{p_j} \cdot s_j,
\end{equation}
where the right hand side is written as additive monoids.

\smallskip
\noindent
{\bf Step 4. Extending the log map to $L \setminus \{q_1, q_2\}$.}
To construct an $f$ over $C$, it remains to construct $f|_{L} \colon L \to X$ which is compatible with \eqref{eq:log-map-away-from-linker} at the two nodes $q_1, q_2$. More precisely, as $\underline{f}(\underline{L}) = x$, $\underline{f}^*\cM_{X}|_{L}$ is generated by the element $\delta$. To define $f|_{L}$, it suffices to find $f^{\flat}(\delta) \in \cM_{C}|_{L}$ whose fibers at the two nodes agree with \eqref{eq:gluing-node-map}. 

We wish to construct the dotted arrows in the following diagram
\[
\xymatrix{
\cT  \ar@{^{(}->}[rr] \ar[d] && \underline{f}^*(\cM_{X})|_{\underline{L}} \ar@{-->}[rr]^{f^{\flat}|_{L}} \ar[d] && \cM_{C}|_{\underline{L}} \ar[d] \\
\{\bar{\delta}\} \ar@{^{(}->}[rr] && \underline{f}^*(\ocM_{X})|_{\underline{L}} \ar@{-->}[rr]^{\bar{f}^{\flat}|_{L}} && \ocM_{C}|_{\underline{L}}
}
\]
where the middle and right vertical arrows are the quotient by $\cO^{\times}_{\underline{L}}$, and the left square is Cartesian. The above discussion implies that $\underline{f}^*(\ocM_{X})|_{\underline{L}} = \mathbb{N}_{\underline{L}}$ is the constant sheaf with the generator $\bar{\delta}$ which is the image of $\delta$. Furthermore $\cT$ is a trivial $\cO^{\times}$-torsor. 

We first define $\bar{f}^{\flat}|_{L}$ by specifying the element $\bar{f}^{\flat}|_{L}(\bar\delta) \in \ocM_{C}|_{\underline{L}}$. Recall that $t_{j},s_{j}$ denote local coordinates at the node $q_{j}$ of $L,Z_{j}$ respectively.  Denote by $\bar{t}_j, \bar{s}_j$ the images of $t_j, s_{j}$ in the characteristic sheaf respectively. Let $t$  be the local coordinate of $p$ which lifts uniquely to a local section of $\cM_{C}|_{\underline{L}}$, again denoted by $t$. Let $\bar{t}$ be the image of $t$ in the characteristic sheaf. 

Note that $\ocM_{C}|_{\underline{L}}$ is constructible with respect to the following strata
\[
q_1, \qquad q_2, \qquad p, \qquad \underline{L}\setminus \{q_1, q_2, p\},
\]
along which $\ocM_{C}$ becomes sheaves of constant monoids. Let $e$ be the generator of $\overline{\cN}$.  We define the image $\bar{f}^{\flat}|_{L}(\bar\delta)$ along each stratum as follows
\begin{equation}\label{eq:stratawais-delta-image}
\begin{aligned}
&\bar{f}^{\flat}|_{L}(\bar\delta)|_{q_j} := c_{p_j} \cdot \bar{s}_j, \qquad \mbox{for } j=1,2; \\
&\bar{f}^{\flat}|_{L}(\bar\delta)|_{p}  := \ell \cdot e + (c_{p_1} + c_{p_2}) \cdot \bar{t};\\
&\bar{f}^{\flat}|_{L}(\bar\delta)|_{\underline{L}\setminus \{q_1, q_2, p\},} := \ell \cdot e.
\end{aligned}
\end{equation}
We observe that the above assignments glue to an element $\bar{f}^{\flat}|_{L}(\bar\delta) \in \ocM_{C}|_{\underline{L}}$. 
First note that the element $(\ell \cdot e + (c_{p_1} + c_{p_2}) \cdot \bar{t})$ generizes to  $\ell \cdot e$ since the element $\bar{t}$ is $0$ when restricted away from $p$. 

To verify that the above construction glues across $q_1$ and  $q_2$,  it suffices to show that $\bar{f}^{\flat}|_{L}(\bar\delta)$ is well-defined in the constructible sheaf of groups $\ocM^{gp}_{C}|_{\underline{L}}$, as all the monoids involved are fine and saturated. By the construction in Step 2, we have 
\[
\ell/c_{p_j} \cdot e =  \bar{s}_j + \bar{t} _j,
\]
hence 
\begin{equation}\label{eq:gluing-node}
c_{p_j} \cdot \bar{s}_j = c_{p_j} \cdot ( \ell/c_{p_j} \cdot  e - \bar{t}_j) = \ell \cdot e - c_{p_j} \cdot \bar{t}_j \qquad \mbox{in } \ocM^{gp}_{C,q_j}. 
\end{equation}
Note that $\bar{t}_j$ is $0$ when away from $q_j$. This shows that $\bar{f}^{\flat}|_{L}(\bar\delta) \in \ocM_{C}|_{\underline{L}}$ is a well-defined section on $\underline{L}$, hence a well-defined morphism on the characteristic level $\bar{f}^{\flat}|_{L} \colon \underline{f}^*\ocM_{X} \to \ocM_{C}|_{L}$.  

Consider the $\cO^{\times}$-torsor $\cT' := \cM_{C}|_{\underline{L}} \times_{\ocM_{C}|_{\underline{L}}}\{\bar{f}^{\flat}|_{L}(\bar\delta)\}$. 
To lift $\bar{f}^{\flat}|_{L} $ to a morphism of log structures $f^{\flat}|_{L} \colon \underline{f}^*\cM_{X}|_{\underline{L}} \to \cM_{C}|_{\underline{L}} $, it suffices to find an isomorphism of $\cO^{\times}$-torsors $\cT \stackrel{\cong}{\longrightarrow} \cT'$.  
Recall that $\cT$ is a trivial torsor. On the other hand,  \eqref{eq:stratawais-delta-image} and \eqref{eq:gluing-node} implies that $\cT'$ consists of sections with poles of order $c_{p_j}$ at $q_j$ and zeros of order $c_{p_1} + c_{p_2}$ at $p$, and no other poles and zeros. In particular,  $\cT'$ is also trivial. Thus, any isomorphism $\cT \cong \cT'$ defines the morphism $f^{\flat}$, hence a log map $f \colon C \to X$ over $S$ as needed.

\subsection{Smoothing of the log gluing}

Let $f \colon C \to X$ over $S$ be constructed as in Section \ref{sss:gluing-along-stratum}. 
We further impose point constraints $f_{P_j}$ for $f_j$ as in \S \ref{ss:point-constraints-deformation}. This leads to point constraints $f_P = f_{P_1} \cup f_{P_2}$ for $f$. 
 We next study the deformation of $f$ with the constraints $f_P$, and prove the following sufficient condition for smoothing $f$ (and in particular smoothing the two nodes $q_1, q_2$).  
 
 \begin{prop}
\label{prop:smoothing lemma}
Notation as above.  Suppose that $c_{p_1}, c_{p_2}, c_p = (c_{p_1} + c_{p_2})$ are coprime to $\ch\mathbf{k}$, the restrictions $H^0(N_{f_j, f_{P_j}}) \to H^0(N_{f_j, f_{P_j}}|_{p_j}) \cong T_{\Delta}|_{x}$ are surjective (see Lemma \ref{lemm:N-at-marking}), and $H^1(N_{f_j, f_{P_j}}) = 0$. Then we have $H^1(N_{f, f_{P}}) = 0$. In particular, a general deformation of $f$ with the point constraints $f_P$ is non-degenerate. 
\end{prop}

 \begin{proof} 
 Note that $H^1(N_{f, f_{P}}) = 0$ implies $f$ is unobstructed relative to $\mathbf{Log}$, see \S \ref{ss:point-constraints-deformation}. Thus, a general deformation of $f$ is a log map over a point with the trivial log structure, i.e. non-degenerate. Thus the last sentence of the statement follows from the earlier claims.
 
For the rest of this statement,  we will assume that $f_P = \emptyset$, hence $f_{P_1} = f_{P_2} = \emptyset$. The case with point constraints is identical but with slightly more complicated notations. 
 
 Consider the following distinguished triangle over $C$
\begin{equation}\label{eq:gluing-cotangent-complex}
f^*\Omega_{X} \stackrel{f^*}{\longrightarrow} \Omega_{C/S} \longrightarrow \mathbb{L}_f \stackrel{[1]}{\longrightarrow}.
\end{equation}
Taking duals and rotating the complex, we obtain a distinguished triangle as in \eqref{eq:deformation-complex}
\begin{equation}\label{eq:dual-gluing-cotangent-complex}
T_{C/S} \stackrel{df}{\longrightarrow} f^*T_X \longrightarrow \mathbb{L}_f^{\vee}[1] \stackrel{[1]}{\longrightarrow}
\end{equation}
with $N_f = \mathbb{L}_f^{\vee}[1]$. 
We would like to compute $H^1(N_f)$.

Similarly we define $\mathbb{L}_{f_j}$ and $N_{f_{j}} = \mathbb{L}^{\vee}_{f_j}[1]$ over $C_j$ for $j=1,2$. 
It follows from the construction in Section \ref{sss:gluing-along-stratum} that 
$
(\mathbb{L}_f)|_{C_j} \cong \mathbb{L}_{f_j},
$
hence $N_{f}|_{C_j} \cong N_{f_{j}}$ for $j=1,2$.

To compute the cohomology of \eqref{eq:dual-gluing-cotangent-complex}, consider the partial normalization sequence 
\[
\pi \colon \tilde{C} := C_1 \sqcup L \sqcup C_2 \to C,
\]
hence a short exact sequence over $C$:
\[
0 \longrightarrow \cO_C \longrightarrow \pi_{*} \cO_{\tilde{C}} \longrightarrow \cO_{q_1}\oplus \cO_{q_2} \longrightarrow 0.
\]
Applying $\otimes N_f$ to the short exact sequence, and taking the corresponding long exact sequence, we have 
\[
\begin{aligned}
0 &\to H^0(N_f) \to & H^0(N_{f_1})\oplus H^0(N_{f}|_{L}) \oplus H^0(N_{f_2})  \to & H^0(N_{f}|_{q_1})\oplus H^0(N_{f}|_{q_2}) & \\
 & \to H^1(N_f) \to & H^1(N_{f_1})\oplus H^1(N_{f}|_{L}) \oplus H^1(N_{f_2})  \to & H^1(N_{f}|_{q_1})\oplus H^1(N_{f}|_{q_2}) & \to \cdots
\end{aligned}
\]

To prove $H^1(N_{f}) = 0$,
we compute $N_{f}|_{L}$ and $N_{f}|_{q_j}$ and their cohomologies as follows.
Consider $(\mathbb{L}_f)|_{L}$ which fits in the distinguished triangle
\[
f^*\Omega_{X}|_{L} \stackrel{f^*|_{L}}{\longrightarrow} \Omega_{C/S}|_{L} \longrightarrow \mathbb{L}_f|_{L} \stackrel{[1]}{\longrightarrow}.
\]
Since $f$ contracts the component $L$ to the point $x \in X$, we have $f^*\Omega_X|_{L} \cong \cO^n_X$ for $n = \dim X$ since $X$ is log smooth. Select local coordinates $y_1, \cdots, y_{n-1}, y_n$ in a neighborhood $U$ of $x$ such that $\Delta\cap U$ is defined by $y_n=0$. Then 
\[
dy_1, ~ \cdots,~ dy_{n-1},~ d\log y_n
\]
form a basis of $\Omega_X$ locally around $x$. By abuse of notations, we view them as a basis of $f^*\Omega_X|_{L}$. On the other hand, we observe that
\[
\Omega_{C/S}|_{L} = \Omega_{\underline{L}}(q_1 + q_2 + p) \cong \cO_{L}(1)
\]
Since $L$ is contracted, we have $f^*|_{L}(dy_k) = 0$ for $k = 1, \cdots n-1$. 

Next we observe that $f^*_{L}(d\log y_n) \in H^0(\Omega_{C/S}|_{L})$ is a non-trivial section that vanishes at a point away from $q_1, q_2$ and $p$. Indeed, since the contact orders $c_{p_1}, c_{p_2}$, and $(c_{p_1} + c_{p_2})$ are all prime to the characteristic of the ground field Equation \eqref{eq:gluing-node} on the characteristic level  implies that locally around $q_j$ we have 
\[
f^*|_{L}(d\log y_n) = u_j \cdot (-c_{p_j}) d\log t_j
\]
for some locally invertible function $u_j$, where $t_j$ is the local coordinate of $L$ around $q_j$ for $j=1,2$. And similarly around $p$ the second equation of \eqref{eq:stratawais-delta-image} implies that 
\[
f^*|_{L}(d\log y_n) = u_p \cdot (c_{p_1} + c_{p_2}) d\log t
\]
for some locally invertible function $u_p$, where $t$ is the local coordinate of $L$ around $p$. As $d\log t_1, d\log t_2$, and $d\log t$ are local generators of $\Omega_{C/S}|_{L}$ at $q_1$, $q_2$ and $p$ respectively, the section $f^*|_{L}(d\log y_n)$ does not vanish at the three special points $q_1, q_2$, and $p$. Thus the restriction of $f^*|_{L}$ to the component $\cO_{L}\cdot d\log y_n \subset f^*\Omega_X|_{L}$ is an injection $\cO_{L} \to \cO_{L}(1) \cong \Omega_{C/S}|_{L}$.
By degree considerations, the section $f^*|_{L}(d\log y_n)$ vanishes at a point other than $q_1, q_2$, and $p$, say $p'$. 

Taking a dual, the above discussion implies that $N_{f}|_{L}$ is the cone of the following complex
\[
\cO_{L}(-1) \stackrel{df}{\longrightarrow} \cO_{L}^{\oplus n}.
\]
In particular, we have $N_{f}|_{L} \cong \cO_{L}^{(n-1)}\oplus \cO_{p'}$. 
Noting that $y_1, \cdots, y_{n-1}$ form local coordinates of $\Delta$ around $x$, by further restricting to $q_j$ we have 
\[
N_{f}|_{q_j} \cong T_{\Delta}|_{x} \cong \cO_{q_j}^{\oplus (n-1)}
\]
for $j=1,2$.  Thus we have the vanishing 
\[
H^1(N_{f}|_{q_1}) = 0, \qquad  H^1(N_{f}|_{q_2}) = 0, \qquad H^1(N_{f}|_{L}) = 0. 
\]

Finally, the statement follows from the above calculations and the long exact sequence, noting that $N_{f_j}|_{p_j} \cong N_{f}|_{q_j}$. 
\end{proof}

\subsection{Splitting contact orders}\label{ss:contact-splitting}

Consider  a log smooth variety $X$ with strict normal crossings boundary $\Delta = \cup_{i \in I} \Delta_i$ as above. 
Consider a non-zero contact order $\mathbf{c}_{\star} = (c_{\star,i})$ 
at a marking, and denote by $I_{{\star}} = \{ i \in I \ | \ c_{\star,i} \neq 0\}$ the set of non-vanishing components of $\mathbf{c}_{\star}$. The {\em total splitting} of $\mathbf{c}_{\star}$ is the collection of non-zero contact orders $\operatorname{split}(\mathbf{c}_{\star}) := \{\mathbf{u}_{k} = (u_{k,i})\}_{k \in I_{\star}}$ where $u_{k,i} = 0$ if $i \neq k$, and $u_{k,k} = c_{\star,k}$.

\begin{prop}
\label{prop:splittingcontactorder}
Suppose that $f \colon C \to X$ is a non-degenerate, rational log map with contact orders $\varsigma = \{\mathbf{c}_k = (c_{k,i})\}_{k}$ such  that  $\ch \mathbf{k} \nmid \mathbf{c}_k$ for all $k$, and point constraints $f_P$ at all non-contact markings. Let $\mathbf{c}_{\star} \in \varsigma$ be a non-zero contact order  at the marked point $p_{\star}$ satisfying $\ch\mathbf{k} \nmid c_{\star,i}$ for all $i$ such that $c_{\star, i} \neq 0$.
Assume $H^{1}(N_{f,f_P}) = 0$. Then there is a non-degenerate, rational log map $\widetilde{f} \colon \widetilde{C} \to X$ with contact orders $\varsigma' = \{ \mathbf{c}_k \ | k \neq \star \} \cup \operatorname{split}(\mathbf{c}_{\star})$ and curve class $f_*[C]$, satisfying the same point constraints $f_P$. 
\end{prop}

\begin{rema}
\label{rema:splittingcontact_freeness}
In the proof, we will construct a degenerate log map $\widetilde{f}_0 \colon \widetilde{C}_0 \to X$ over $S$ with contact orders $\varsigma'$ and point constraints $f_P$ such that $H^{1}(N_{\widetilde{f}_0,f_P}) = 0$ and $f \colon C \to X$ appears as one of its components. Then $\widetilde{f}$ is constructed as a general smoothing of $\widetilde{f}_0$. In some sense, $\widetilde{f}$ can be viewed as a ``deformation'' of $f$ that splits the contact order $\mathbf{c}_{\star}$. 
By upper semicontinuity of cohomologies of complexes of coherent sheaves which are flat over the base (see for example \cite[Proposition 6.4]{Hartmann12}), we conclude that a general smoothing of $\widetilde{f}_0$ satisfies $H^1(N_{\widetilde{f},f_P}) = 0$.
\end{rema}

\begin{proof}
We will assume in the proof that $f_{P} = \emptyset$. The general case is similar but with more complicated notation and is left to the readers. 

Denote by $I_{\star} = \{ i \in I \ | \ c_{\star,i} \neq 0\}$.  
If $|I_{\star}| = 1$, then there is nothing to prove. So we assume $d := |I_{\star}| \geq 2$. We split the proof into several steps. 

{\bf Step 1: Construct a degenerate domain curve.} 
Consider $\underline{R} \cong \mathbb{P}^1$ with a choice of distinct points $q_{\star}, q_1, q_{2}, \cdots, q_{|I_{\star}|} \in \underline{R}$. 
Let $\underline{\widetilde{C}}_0 = \underline{C} \cup \underline{R}$ obtained by identifying $p_{\star} \in \underline{C}$ with $q_{\star} \in \underline{R}$. This pre-stable curve $\underline{\widetilde{C}}_0$ has the set of markings $\{p_k\}_{k \neq \star} \cup \{q_{1}, q_2, \cdots, q_{|I_{\star}|}\}$. 
We obtain a log curve $\widetilde{C}_0 \to S$ with the underlying pre-stable curve $\underline{\widetilde{C}}_0$ and its canonical log structure uniquely determined by its underlying structure.

\smallskip

{\bf Step 2: Construct a degenerate log map.}
Consider the stable map $\underline{\widetilde{f}}_0 \colon \underline{\widetilde{C}}_0 \to \underline{X}$ such that $\underline{\widetilde{f}}_0|_{\underline{C}} = \underline{f}$ and $\underline{\widetilde{f}}_0(\underline{R}) = \underline{f}(p_{\star})$. We will lift  $\underline{\widetilde{f}}_0$ to a stable log map $\widetilde{f}_0 \colon \widetilde{C}_0 \to X$ such that its contact orders at $\{q_{1}, q_2, \cdots, q_{|I_{\star}|}\}$ are given by  $\operatorname{split}(\mathbf{c}_{\star}) = \{\mathbf{u}_{1}, \mathbf{u}_{2}, \cdots, \mathbf{u}_{|I_{\star}|}\}$, where $\mathbf{u}_{k}$ mean the non-trivial contact order with respect to $\Delta_{k}$.  

For each $i \in I$, denote by $X_i$ the log scheme associated to the pair $(\underline{X}, \Delta_i)$. If $i \in I_{\star}$, we may apply \cite[Lemma 3.6]{CZ14} along the component $\underline{R}$ to obtain a stable log map $h_i \colon \widetilde{C}_0 \to X^{\dagger}_i$ such that $q_{i}$ has contact order  $\mathbf{c}_{\star,i}$. 
Here $X^{\dagger}_i$ is the specialization to normal cone of $\Delta_i \subset \underline{X}$ with its canonical log structure and a canonical projection $X^{\dagger}_i \to X_i$.  
Now let $\widetilde{f}_{0,i}$ be the stable log map given by the composition $\widetilde{C}_0 \to X^{\dagger}_i \to X_i$. 

If $i \not\in I_{\star}$, since $\underline{\widetilde{f}}_0|_{\underline{C}} = \underline{f}$ is non-degenerate, the underlying structure naturally induces a log map $\widetilde{f}_{0,i} \colon \widetilde{C}_0 \to X_i$. 
Finally, we define 
\[
\widetilde{f}_0 := \prod_{i \in I} \widetilde{f}_{0,i} \colon \widetilde{C}_0 \to X = \prod_{i\in I}X_i
\]
where the product $ \prod_{i\in I}X_i$ is taken over $\underline{X}$. Observe that $\widetilde{f}_0$ has contact orders $\varsigma'$. 

\smallskip
{\bf Step 3: Smoothing of the degenerate log map.} We will show that $H^1(N_{\widetilde{f}_0}) = 0$. Hence $\widetilde{f}_0$ can be deformed to a non-degenerate log map with contact orders $\varsigma'$. We work out the details following the same line of calculation as in Proposition \ref{prop:smoothing lemma}. 

Similar to \eqref{eq:gluing-cotangent-complex} and \eqref{eq:dual-gluing-cotangent-complex}, we have two triangles
\begin{equation}\label{eq:split-contacts-triangle}
\widetilde{f}_0^*\Omega_{X} \stackrel{\widetilde{f}_0^*}{\longrightarrow} \Omega_{\widetilde{C}_0/S} \longrightarrow \mathbb{L}_{\widetilde{f}_0} \stackrel{[1]}{\longrightarrow}, \qquad T_{\widetilde{C}_0/S} \stackrel{d \widetilde{f}_0}{\longrightarrow} \widetilde{f}_0^*T_X \longrightarrow \mathbb{L}_{\widetilde{f}_0}^{\vee}[1] \stackrel{[1]}{\longrightarrow},
\end{equation}
such that $\mathbb{L}_{\widetilde{f}_0}^{\vee}[1] \cong N_{\widetilde{f}_0}$. Using the normalization $\underline{C}\sqcup \underline{R} \to \underline{\widetilde{C}}_{0}$ and arguing as in Proposition \ref{prop:smoothing lemma}, we obtain a long exact sequence
\[
\begin{aligned}
0 &\to H^0(N_{\widetilde{f}_0}) \to & H^0(N_{f})\oplus H^0(N_{\widetilde{f}_0}|_{R})  \to & H^0(N_{\widetilde{f}_0}|_{q_\star}) & \\
 & \to H^1(N_{\widetilde{f}_0}) \to & H^1(N_{f})\oplus H^1(N_{\widetilde{f}_0}|_{R})   \to & H^1(N_{\widetilde{f}_0}|_{q_\star}) & \to \cdots
\end{aligned}
\]
noting that $N_{f} \cong N_{\widetilde{f}_0}|_{C}$. We next compute $N_{\widetilde{f}_0}|_{R}$ and its cohomologies. 

Select local coordinates $y_1, \cdots, y_{n-1}, y_n$ in a neighborhood $U$ of $f(p_{\star})$ such that $\Delta_i\cap U$ is defined by $y_i=0$ for $i \in I_{\star}$. Then we have a basis of $\Omega_X|_{U}$:
\[
\{dy_i \ | \  i \not\in I_{\star} \} \sqcup \{d\log y_i \ | \ i \in I_{\star}\}
\]
By an abuse of notation, we view them as a basis of $\widetilde{f}_0^*\Omega_X|_{R}$. Since $R$ is contracted, we obtain 
\[
\widetilde{f}_0^*\Omega_{X}|_{R} \cong \cO_R^{\oplus n- |I_{\star}|} \oplus \cO_R^{\oplus|I_{\star}|}
\]
given by the above choice of basis.  

We also observe that
$
\Omega_{\widetilde{C}_0/S}|_{R} = \Omega_{\underline{R}}(q_{\star} + \sum_{i \in I_{\star}} q_{i}) \cong \cO_{R}(|I_{\star}| - 1). 
$
Since $R$ is contracted, we have $\widetilde{f}_0^*(dy_i)|_{R} = 0$ for $i \not\in I_{\star}$. Furthermore, since $\ch\mathbf{k} \nmid c_{\star,i}$ for all $i$, we check that for each $i \in I_{\star}$ the section $\widetilde{f}_0^*(d\log y_i)|_{R} \in \Omega_{\widetilde{C}_0/S}|_{R}$ is a local generator at $q_i$, but is not a local generator at $q_j$ for $i \neq j$. Thus $\{ \tilde{f}_0^*(d\log y_i)|_{R} \ | \ i \in I_{\star} \}$ form a basis of $H^0(\Omega_{\tilde{C}_0/S}|_{R})$. In particular, the restriction $\tilde{f}_0^*|_{R}$ is a surjection of vector bundles
\[
\widetilde{f}_0^*\Omega_{X}|_{R} \cong \cO_R^{\oplus n- |I_{\star}|} \oplus \cO_R^{\oplus |I_{\star}|} \to \Omega_{\widetilde{C}_0/S}|_{R} \cong \cO_{R}(|I_{\star}| - 1).
\]
Taking duals, we obtain an exact sequence 
\[
0 \to \cO_{R}(-|I_{\star}| + 1) \to \cO_R^{\oplus n- |I_{\star}|} \oplus \cO_R^{\oplus |I_{\star}|} \to N_{\widetilde{f}_0}|_{R} \to 0
\]
realizing $N_{\widetilde{f}_0}|_{R}$ as a semi-positive vector bundle over $R$. In particular, the restriction morphism $H^0(N_{\widetilde{f}_0}|_{R})  \to  H^0(N_{\widetilde{f}_0}|_{q_\star})$ is surjective, and $H^1(N_{\widetilde{f}_0}|_{R}) = 0$. 

Finally applying the assumption $H^{1}(N_{f}) = 0$ and the long exact sequence, we obtain $H^1(N_{\widetilde{f}_0}) = 0$ as needed. 
\end{proof}

\section{Deformation theory for log sections}
\label{sec:deformation_logsection}

In this section we analyze the sections of a morphism $\pi: \mathcal{X} \to B$ where $\mathcal{X}$ is a log scheme defined by a strict normal crossing divisor and $B$ is a curve.  The main goal is to understand the moduli space and deformation theory of sections satisfying various properties.

\subsection{Stable log sections and their stacks}
Let $\cX$ be a log scheme such that $\underline{\cX}$ is smooth and the boundary divisor $\Delta \subset \underline{\cX}$ is strict normal crossings. A flat, generically smooth, and projective morphism $\pi \colon \cX \to B$ with connected fibers is called a {\em log fibration} if furthermore the restriction $\pi|_{\Delta} \colon \Delta \to B$ is flat and generically relatively strict normal crossings. 
Here we will only consider the case that $B$ is a proper, smooth, genus $g$ curve equipped with the trivial log structure. 

Let $S$ be a log point.
A genus $g$ stable log map $f \colon C \to \cX$ over $S$ is called a {\em stable log map of section type} of $\pi$ if 
 the curve class of the underlying pre-stable map of the composition $C \to \cX \to B$ is $[B]$. 
 It is further called a {\em log section} if $f$ is non-degenerate.

Denote by $\overline{\Sec}_{\log}(\cX/B)$ the stack of stable log maps of section type for the log fibration $\pi \colon \cX \to B$, and $\overline{\Sec}_{\log}(\cX/B)^{+} \subset \overline{\Sec}_{\log}(\cX/B)$ the open and closed substack parametrizing log maps with only contact markings. 
 The open substacks 
 \[
 \Sec_{\log}(\cX/B) \subset \overline{\Sec}_{\log}(\cX/B) \qquad \mbox{and} \qquad \Sec_{\log}(\cX/B)^{+} \subset \overline{\Sec}_{\log}(\cX/B)^{+}
 \] 
 with the trivial log structure are the corresponding moduli stacks of log sections. Note that these open substacks are honest schemes. 
 There is a decomposition into open and closed substacks
 \[
 \overline{\Sec}_{\log}(\cX/B) = \bigsqcup_{(\varsigma, \beta)} \overline{\Sec}_{\varsigma}(\cX/B, \beta)  
  \]
 by further specifying the number of markings, contact orders $\varsigma$, and curve classes $\beta \in N_1(\cX)$
 with $\pi_*\beta = [B]$. We observe that $\overline{\Sec}_{\varsigma}(\cX/B, \beta) = \mathscr{M}_{g, \varsigma}(\cX,\beta)$.  Similarly we have  
 \[
 \overline{\Sec}_{\log}(\cX/B)^{+} = \bigsqcup_{(\varsigma, \beta)} \overline{\Sec}_{\varsigma}(\cX/B, \beta)
\]
is the union over $\varsigma$ without non-contact marking.

\subsection{Deformations of stable log maps of section type}

Let $[f] \in \overline{\Sec}_{\varsigma}(\cX/B, \beta)(S)$ be an object 
over a log scheme $S$. Consider the commutative triangle
\begin{equation}\label{diag:section}
\xymatrix{
 && \cX \ar[d] \\
C \ar[rru]^{f} \ar[rr]^{h} && B
}
\end{equation}
This induces a tautological morphism 
\begin{equation}\label{eq:sec-relative-moduli}
\overline{\Sec}_{\varsigma}(\cX/B,\beta)  \to \mathfrak{M}_{g, |\varsigma|}(B, [B])
\end{equation}
given by $[f \colon C \to \cX] \mapsto [h \colon C \to B]$. Here  $\mathfrak{M}_{g, |\varsigma|}(B, [B])$ is the Artin stack of $|\varsigma|$-marked, genus $g(B)$ pre-stable log maps to the target $B$ with curve class $[B]$. 
It fits in a commutative diagram with strict top arrows:
\begin{equation}\label{eq:forget}
\xymatrix{
\overline{\Sec}_{\varsigma}(\cX/B,\beta) \ar[r] \ar[rd] & \mathbf{Log}_{\mathfrak{M}_{g, |\varsigma|}(B, [B])} \ar[r] \ar[d] & \mathbf{Log} \ar[d] \\
&  \mathfrak{M}_{g, |\varsigma|}(B, [B]) \ar[r] & \Spec\ \mathbf{k}
 }
\end{equation}
where   $\mathbf{Log}_{\bullet}$ is Olsson's log stack parametrizing log structures over a given log stack $\bullet$ as in \S \ref{ss:deformation}.

Consider an object $[f] \in \overline{\Sec}_{\varsigma}(\cX/B,\beta)(S)$ as in \eqref{diag:section} over a log point $S$.  The normal complexes $N_{f}, N_{h}$ as in \S \ref{ss:deformation} fit in the following commutative diagram of distinguished triangles:
\begin{equation}\label{diag:deformation-octahedron}
\xymatrix{
T_{C/S} \ar[rrd]^{d f} \ar[rdd]_{d h} && && \\
&& f^*T_{\cX}  \ar[rrd] \ar[ld] && \\
& h^*T_{B} \ar[ldd] \ar[rd] &&& N_{f} \ar[lld] \\
&& N_h \ar[lld] &&  \\
f^*T_{\cX/B}[1] && &&
}
\end{equation}

In particular, we have the triangle
\begin{equation}\label{eq:deformation-triangle}
f^*T_{\cX/B} \longrightarrow N_f \longrightarrow N_h \stackrel{[1]}{\longrightarrow}
\end{equation}

By \eqref{eq:deformation-complex}, we see that the cohomologies of $N_{f}$ and $N_{h}$ control the deformations of $f$ and $h$ relative to $\mathbf{Log}$ respectively. 

Next, consider the deformations of \eqref{diag:section} relative to $\mathbf{Log}_{\mathfrak{M}_{g, |\varsigma|}(B, [B])}$. These are deformations of $[f]$ while fixing $([h], S)$). In this case the first-order deformations and obstructions of \eqref{diag:section} relative to $\mathbf{Log}_{\mathfrak{M}_{g, |\varsigma|}(B, [B])}$ are given by $H^0(f^*T_{\cX/B})$ and $H^1(f^*T_{\cX/B})$ respectively.

Under certain conditions the objects in the derived category defined by Equation \eqref{diag:deformation-octahedron} are represented by sheaves.

\begin{prop} \label{prop:normalsheafsequence}
Consider a log section $[f] \in {\Sec}_{\log}(\cX/B)(S)$ where $S$ is a geometric point with the trivial log structure (whose image is thus necessarily a section of $\pi\colon \mathcal{X} \to B$).  Let $\{p_k\}_{k \in P}$ be the collection of non-contact markings and let $\bar{f} \colon \bar{C} \to \cX$ be log map associated to $f$ by removing markings in $\{p_k\}_{k \in P}$, see \S \ref{ss:immersion-deformation}.  We further assume that $\ch\mathbf{k}$ does not divide any of the non-zero contact orders of $f$.  Then
\begin{enumerate}
 \item $\bar{f}$ is a log immersion.
 \item $N_{f} \cong N_{\bar{f}}\oplus\bigoplus_{k\in P}N_{p_k/\underline{C}}$ and $N_h \cong \oplus_{k=1}^{|\varsigma|} N_{p_k/\underline{C}}$ are sheaves, where $\{p_1, \cdots, p_{|\varsigma|}\}$ is the set of all markings and $N_{p_i/C} \cong \cO_{p_i}$ is the normal bundle of the underlying marking.
 \item Further suppose that $f(C)$ is contained in the locus where $\pi: \mathcal{X} \to B$ is log smooth (where $B$ is equipped with the trivial log structure).  Then we have an exact sequence
\begin{equation}\label{eq:deformation-sequence}
0 \longrightarrow f^*T_{\cX/B} \longrightarrow N_f \longrightarrow N_h = \oplus_{i=1}^{|\varsigma|} N_{p_i/C} \longrightarrow 0
\end{equation}
where $N_{p_i/C} \cong \cO_{p_i}$ is the normal bundle of the underlying marking.
\end{enumerate}

\end{prop}

\begin{proof}

Noting that the underlying morphism of $f$ is a closed embedding, (1) follows from Definition \ref{defi:log-immersion}.  (2) is a consequence of Corollary \ref{cor:torsion-normal-complex}.

(3) 
We have the following commutative diagram
\[
\xymatrix{
0 \ar[r] &T_{C} \ar[r] \ar[d]_{=} &  f^*T_{\cX} \ar[r] \ar[d] & N_f \ar[r] \ar[d] & 0 \\
0 \ar[r] & T_{C} \ar[r] & T_{B} \ar[r] & \oplus_i N_{p_i/C} \ar[r] & 0
}
\]
Since we are assuming $C$ is contained in the log smooth locus of $\pi$, the map $f^{*}T_{\cX} \to T_{B}$ is surjective.  The kernel-cokernel sequence shows that the rightmost arrow is also surjective and has kernel isomorphic to $T_{\mathcal{X}/B}$. 
\end{proof}

\subsection{Log sections through given points}
Consider a finite set of points 
\[
\{q_k\}_{k \in P} \subset \mathcal{X}^\circ = \mathcal X \setminus \mathrm{Supp}(\Delta)
\]
such that their images $p_{k} = \pi(q_{k})$ are distinct. 
Denote by 
\[
\Sec_{\log}(\cX/B, \{p_k\}_{k \in P}) \subset \Sec_{\log}(\cX/B)
\] 
the moduli of log sections with a specified subset of additional markings with images $\{p_k\}_{k \in P}$ in $B$.  By a mild abuse of notation, we will use $p_{k}$ to denote the corresponding marking of the domain curve $C$ (as well as its image in $B$). Denote by 
\[
\Sec_{\log}(\mathcal{X}/B,\{q_k\}_{k \in P}) \subset \Sec_{\log}(\cX/B, \{p_k\}_{k \in P})
\] 
the closed substack parametrizing log sections such that the image of the marking $p_k$ is $q_k$ for all $k$.
This is the moduli space of log sections through $\{q_{k}\}_{k \in P} \subset \cX$. Denote by 
\begin{align*}
\Sec_{\log}(\mathcal{X}/B,\{p_k\}_{k \in P})^{+} & \subset \Sec_{\log}(\mathcal{X}/B,\{p_k\}_{k \in P}\}) \\
\Sec_{\log}(\mathcal{X}/B,\{q_k\}_{k \in P})^{+} & \subset \Sec_{\log}(\mathcal{X}/B,\{q_k\}_{k \in P})
\end{align*}
the open and closed substacks such that the markings with zero contact orders are exactly $\{p_{k}\}_{k \in P}$.

Let $[f] \in \Sec_{\log}(\mathcal{X}/B,\{q_k\}_{k \in P})^{+}(S)$ be a log section over a geometric log point $S$ with the trivial log structure, and $\bar{f} \colon \bar{C} \to \cX$ be log map associated to $f$ by removing markings in $\{p_k\}_{k \in P}$. 
As in \eqref{eq:twisted-normal-complex}, we define the {\em twisted normal complex}:
\[
N_{f, \{q_k\}} := N_{\bar{f}}\left(-\sum_{k \in P} q_k \right).
\]
The cohomology groups $H^0(N_{f, \{q_k\}})$ and $H^1(N_{f, \{q_k\}})$ control the first-order infinitesimal deformations and obstructions of $[f]$ inside of $\Sec_{\log}(\mathcal{X}/B,\{q_k\}_{k \in P})$ relative to $\mathbf{Log}$, see \S \ref{ss:point-constraints-deformation}. 
Further assuming that $\ch \mathbf k$ does not divide any of the non-zero contact orders of $f$, then $N_{\bar{f}}$ and $N_{f, \{q_k\}}$ are vector bundles by Proposition \ref{prop:normalsheafsequence}.

\begin{lemm}
\label{lemm:goingthroughgeneralpoints}
Fix points $\{q_k\}_{k \in P} \subset \mathcal{X}^\circ$ with distinct images $p_{k} = \pi(q_{k})$ in $B$.

Suppose that $M$ is an irreducible substack 
of $\Sec_{log}(\mathcal{X}/B,\{q_{k}\}_{k\in P})^{+}$ 
parametrizing a separable dominant family of log sections such that all
of the non-zero contact orders are not divisible by $\ch \mathbf k$.  Then for a general member $[f: C \to \mathcal{X}] \in M(\Spec \ \mathbf{k})$, the vector bundle $N_{f, \{q_k\}}$ is generically globally generated.

Conversely, suppose we fix a log section $f: C \to \mathcal{X}$ through $\{q_k\}_{k \in P}$ such that $N_{f, \{q_k\}}$ is generically globally generated and $H^{1}(C,N_{f, \{q_k\}})= 0$.  
Then there is a unique irreducible component $M \subset \Sec_{log}(\mathcal{X}/B,\{q_{k}\}_{k\in P})^{+}$ 
containing $[f]$. Furthermore $M$ parametrizes a separable dominant family of log sections.
\end{lemm}
\begin{proof}
We first prove the first statement.  Let $\pi : \mathcal U \to M$ be the universal family 
with the evaluation map $ev : \mathcal U \to \cX$. For a point $(f : C \to \mathcal X, p) \in \mathcal U$ with $p \in C$ non-marked, the tangent space is given by $TM_{f} \oplus T_C|_p$. Since $ev$ is dominant and separable, for a general choice of $f$ and $p$ we have a surjection
\[
TM_{f} \oplus T_C|_p \subset H^0\left(C, N_{f, \{q_k\}}\right) \oplus T_C|_p \to f^*T_{\mathcal X}|_p.
\]
This implies that 
\[
H^0\left(C, N_{f, \{q_k\}}\right)\to N_{f, \{q_k\}}|_p
\]
is surjective. Thus $N_{f, \{q_k\}}$ is generically globally generated.

Next we prove the second statement; the proof is almost backward. For a general non-marked point $p \in C$, 
\[
H^0\left(C, N_{f, \{q_k\}}\right)\to N_{f, \{q_k\}}|_p
\]
is surjective proving that
\[
H^0\left(C, N_{f, \{q_k\}}\right) \oplus T_C|_p \to f^*T_{\mathcal X}|_p.
\]
is surjective. Furthermore, $H^{1}(C,N_{f, \{q_k\}})= 0$ implies that $\mathcal U$ is smooth around the point $p$. 
Thus $ev : \mathcal U \to \mathcal X$ is dominant and separable, proving the claim.
\end{proof}

\section{Free log curves}
\label{sec:freelogcurve}

A key tool for understanding rational curves on projective varieties is the notion of a (very) free curve.  The study of free curves is closely tied to the geometric notions of uniruledness and rational connectedness.  In this section we quickly describe the analogous theory in the log setting.  Similar proposals have been put forward by e.g.~\cite{KM99,Campana07,Campana10,CZ15, CZ17, CZ18}.

\subsection{Uniruledness and connectedness by rational log curves}

We introduce the notation of (separable) uniruledness and rational connectedness of $X$ by log curves.

\begin{defi}\label{def:uniruled-connected}
Let $X$ be a log smooth variety with strict normal crossings boundary and let $\varsigma$ be a collection of non-zero contact orders.  
\begin{enumerate}
\item We say $X$ is (separably) $\varsigma$-uniruled if there is a family of genus zero non-degenerate stable log maps $\pi : \mathcal U \to W, ev \colon \mathcal U \to X$ with contact orders $\varsigma$ such that $\underline{\mathcal U} = W \times \mathbb{P}^1$, $\dim W = \dim X - 1$, and $ev$ is dominant (and separable). 

\item We say $X$ is (separably) $\varsigma$-rationally connected if there is a family of genus zero non-degenerate stable log maps $\pi : \mathcal U \to W, ev \colon \mathcal U \to X$ with contact orders $\varsigma$ such that  $\underline{\mathcal U} = W \times \mathbb{P}^1$ and $ev^{(2)} \colon \mathcal U \times_{W} \mathcal U \to X^2$ is dominant (and separable).
\end{enumerate}
\end{defi}

\begin{defi}\label{def:free-curve}
Notations as in Definition \ref{def:uniruled-connected}, let $f: C \to X$ be a non-degenerate rational log curve over a geometric point with contact orders $\varsigma$. We say that $f$ is {\em $\varsigma$-free} (resp. {\em $\varsigma$-very free}) if $H^1(N_f(-1)) = 0$ (resp. $H^1(N_f(-2)) = 0$).  
\end{defi}

\begin{coro}\label{cor:easy-free}
Let $f: C \to X$ be a non-degenerate rational log curve over a geometric point with contact orders $\varsigma$. If $f^*T_X$ is nef (resp. ample), then $f$ is free (resp.~very free).  
\end{coro}
\begin{proof}
To show that $f$ is $\varsigma$-free (resp.~$\varsigma$-very free), one may first apply $\otimes \cO_{\mathbb{P}^1}(-1)$ (resp. $\otimes \cO_{\mathbb{P}^1}(-2)$) to \eqref{eq:deformation-complex} and then take the corresponding long exact sequence. The statement then follows from Definition \ref{def:free-curve}. 
\end{proof}

\begin{prop}\label{prop:uniruled=free}
Let $X$ be a log smooth variety with strict normal crossings boundary and let $\varsigma$ be a collection of non-zero contact orders not divisible by $\ch\mathbf{k}$. Then $X$ is separably $\varsigma$-uniruled (resp.~separably $\varsigma$-rationally connected) iff it admits a $\varsigma$-free (resp.~$\varsigma$-very free) rational log curve.
\end{prop}

We postpone the proof to Section \ref{ss:equivalence}.

\subsection{Relatively free log sections}
Recall that $B$ is of genus $g$. 
The above definition
suggests the following notions of freeness for log sections:

\begin{defi}\label{def:free-section}
Let $f: C \to \mathcal{X}$ be a log section with only contact markings. 
Assume that every contact order is not divisible by $\ch \mathbf k$.
It is called {\em relatively generically free} if $H^{1}(C, N_{f}) = 0$ and $N_{f}$ is generically globally generated.  If furthermore $N_{f}$ is globally generated, it is called {\em relatively free}. 
Note that it follows from Proposition~\ref{prop:normalsheafsequence} that the normal complex $N_f$ is actually a sheaf on $C$. 

We call such an $f$ {\em relatively HN-free} (resp. relatively HN-very free) if $\mu^{min}(N_f) \geq 2g$ (resp. $\mu^{min}(N_f) \geq 2g+1$) where $\mu^{min}(\mathcal F)$ is the minimal slope of a torsion free sheaf $\mathcal F$.

\end{defi}

\begin{lemm}
Let $f: C \to \mathcal{X}$ be a log section. Then we have the following implications:
\[
\text{$f$ is relatively HN-free} \implies \text{$f$ is relatively free} \implies \text{$f$ is relatively generically free.}
\]

In case $g = 0$, the three notions of freeness are all equivalent.  
Thus, we say $f$ is {\em relatively free} if one of the three equivalent conditions is satisfied. Furthermore, when $g = 0$ we say $f$ is {\em relatively very free} if it is relatively HN-very free.
\end{lemm}

\begin{proof}
The second implication is trivial. The first implication follows from Riemann-Roch and Serre duality.  More precisely, it is an immediate consequence of \cite[Corollary 2.8]{LRT23}; the cited paper works over $\mathbb{C}$ but this argument is valid in arbitrary characteristic.  
\end{proof}

\begin{lemm}
\label{lemm:generalpoints}
Let $f : C \to \mathcal X$ be a log section with only contact markings and assume that all of the contact orders of $f$ are not divisible by $\ch \mathbf k$. 
Let $M \subset \Sec_{log}(\mathcal X/B)$ be a component containing $f$.
\begin{enumerate}
\item Suppose that $f$ is relatively HN-free and let $b = \mu^{min}(N_f)$. Then deformations of $C$ parametrized by $M$ go through $\lfloor b \rfloor -2g + 1$ general points;
\item Conversely suppose that our ground field has characteristic $0$ and deformations of $f$ parametrized by $M$ go through $2g+1$ general points.  Then a general deformation of $f$ in $M$ is relatively HN-free.
\end{enumerate}
\end{lemm}

\begin{proof}
(1) Let $p_1, \cdots, p_m$ be general points on $B$. Then as long as $m \leq \lfloor b \rfloor  -2g$, the vector bundle $N_f(-p_1 - \cdots -p_m)$ is globally generated with vanishing $H^1$ by \cite[Corollary 2.8]{LRT23}. Thus our assertion follows from Lemma~\ref{lemm:goingthroughgeneralpoints}.

(2) We fix a deformation $f' : C'\to \mathcal X$ of $f$ going through a set of $2g(B)$ general points on $\mathcal X$. Let $p_1, \cdots, p_{2g}$ be the corresponding points on $C'$. Then it follows from  Lemma~\ref{lemm:goingthroughgeneralpoints} that $N_{f'}(-p_1 - \cdots -p_{2g})$ is generically globally generated. \cite[Lemma 2.6]{LRT23} shows that $\mu^{min}(N_{f'}(-p_1 - \cdots -p_{2g})) \geq 0$. Thus it follows that $\mu^{min}(N_{f'})\geq 2g$.
\end{proof}

\subsection{(Very) free log curves via (very) free log sections}\label{ss:equivalence}
Next, we relate the notions of free log curves
and log sections in Definitions \ref{def:free-curve} and \ref{def:free-section}. 

Given a non-degenerate stable log map $f \colon C \to X$ with contact orders $\varsigma$ over a geometric point $S$, it naturally induces a log section of a trivial family:
\[
\xymatrix{
 && X \times \underline{C} \ar[d]^{\pi} \\
 C \ar[rr]_{u} \ar[rru]^{\rho_f := f\times u} && \underline{C}
}
\]
where $u$ is the morphism forgetting the log structure.  We say the induced log section is {\em non-trivial} if the composition $C \stackrel{\rho}{\longrightarrow} X \times \underline{C}  \stackrel{\pi}{\longrightarrow}  X$ is not a contraction of the curve. Note that $\rho_{f}$ is non-trivial iff the image $f(C)$ is not a point. In particular $\rho_f$ is non-trivial if $\varsigma$ is not entirely zero. 

\begin{lemm}\label{lem:free-curve=free-section}
Notation as above. 
Suppose $C$ is rational and $\varsigma$ consists of non-zero contact orders which are not divisible by $\mathrm{char}\mathbf k$. Then $f$ is free (resp. very free) iff $\rho_{f}$ is relatively HN-free (resp. relatively HN-very free).  
\end{lemm}
\begin{proof}
Consider the following commutative diagram
\[
\xymatrix{
T_{C} \ar[rrd] \ar[rdd] && && \\
&& \rho_{f}^*T_{\cX}  \ar[rrd] \ar[ld] && \\
& f^*T_{X} \ar[ldd] \ar[rd] &&& N_{\rho_f} \ar[lld] \\
&& N_f \ar[lld] &&  \\
T_{\underline{C}}[1] && &&
}
\]
This provides a distinct triangle
\[
T_{\underline{C}} \lra N_{\rho_f} \lra N_f \stackrel{[1]}{\lra} T_{\underline{C}}[1].
\]
Since $\underline{C} \cong \mathbb{P}^1$, by applying $\otimes \cO_{\underline{C}}(-m)$ to the above distinct triangle for $m=1, 2$ and taking the associated long exact sequence, we obtain an exact sequence
\[
0 \lra H^1(N_{\rho_f}(-m)) \lra H^1(N_{f}(-m)) \lra 0
\]

Suppose $f$ is free (resp.~very free), i.e. $H^1(N_{f}(-1)) = 0$ (resp. $H^1(N_{f}(-2)) = 0$) by Definition \ref{def:free-curve}. By the above exact sequence, this is equivalent to $H^1(N_{\rho_f}(-1)) = 0$ (resp. $H^1(N_{\rho_f}(-2)) = 0$).  Since $C$ is a rational curve, this is equivalent to assuming that $\rho_f$ is relatively HN-free (resp.~relatively HN-very free) as in Definition \ref{def:free-section}. This concludes the proof. 
\end{proof}

\begin{lemm}
\label{lemm:varsigma-rat_connected}
Let $X$ be a log smooth projective variety with strict normal crossings boundary, and let $\pi \colon X \times \mathbb{P}^1 \to \mathbb{P}^1$ be the trivial log fibration. Suppose $\varsigma$ consists of non-zero contact orders not divisible by $\ch \mathbf{k}$.
Then we have
\begin{enumerate}
\item $X$ is separably $\varsigma$-uniruled if and only if the log fibration $\pi$ admits a non-trivial, free log section with contact orders $\varsigma$;
\item $X$ is separably $\varsigma$-rationally connected if and only if the log fibration $\pi$ admits a very free log section with contact orders $\varsigma$.
\end{enumerate}
\end{lemm}

\begin{proof}

(1) Suppose that $X$ is separably $\varsigma$-uniruled. Then there is a family $\pi : \mathcal U \to W$ of log rational curves with contact orders $\varsigma$ such that $\underline{\mathcal U} = W \times \mathbb P^1$, $\dim W = \dim X - 1$, and the evaluation map $s : \mathcal U \to X$ is dominant and separable. Let $\phi \colon \mathcal U \to \mathbb{P}^1$ 
be the projection. We thus obtain a family of non-trivial log sections over $W$: 
\begin{equation}\label{eq:map-to-section}
\xymatrix{
 && X \times \mathbb{P}^1 \ar[d]^{\pi} \\
\mathcal{U} \ar[rr]_{\phi} \ar[rru]^{(f,\phi)}&& \mathbb{P}^1 
}
\end{equation}
While this family is not dominant (since $\dim(\mathcal{U}) < \dim(X) + 1$) we may select a torus action to sweep out the total space as follows.
Fix two distinct points $x_1, x_2 \in \mathbb{P}^1$ and an identification $\mathbb{P}^1 \setminus \{x_1, x_2\} \cong \mathbb{G}_m$ with a torus. This defines a multiplication morphism
\[
\mathbf{m} \colon \mathbb{G}_m \times \mathbb{P}^1 \to \mathbb{P}^1
\]
with two fixed points $\{x_1, x_2\}$. This $\mathbf{m}$ is a separable and dominant morphism. Now we extend \eqref{eq:map-to-section} to the following commutative diagram
\[
\xymatrix{
 && X \times \mathbb{G}_m \times \mathbb{P}^1 \ar[d] \ar[rr] && X \times \mathbb{P}^1 \ar[d]^{\pi}  \\
\mathbb{G}_m\times \mathcal{U} \ar[rr] \ar[rru]&& \mathbb{G}_m \times \mathbb{P}^1 \ar[rr]^{\mathbf{m}} && \mathbb{P}^1 
}
\]
The composition of the upper arrows gives a family of non-trivial log sections
\[
\widetilde{ev} \colon \widetilde{\mathcal{U}} := \mathbb{G}_m\times\mathcal{U} \to X\times \mathbb{P}^1
\]
over $\widetilde{W} := \mathbb{G}_m\times W$.  Note that $\widetilde{ev}$ is separable and dominant.  
Let $f : C \to X \times \mathbb P^1$ be a general log section parametrized by $\widetilde{W}$.
Lemma~\ref{lemm:goingthroughgeneralpoints} 
implies that $N_f$ is generically globally generated. Since the underlying curve is a smooth rational curve, we conclude that $f$ is free.

Conversely, suppose that there is a free log section $f : C \to X \times \mathbb P^1$ with contact orders $\varsigma$. 
Let $M$ be the irreducible component of $\Sec_{\varsigma}(X\times \mathbb P^1/\mathbb P^1)$ containing $f$. Let $\pi : \mathcal U \to M$ be the universal family over $M$ with the evaluation map $ev : \mathcal U \to X \times \mathbb P^1$. 
 Lemma~\ref{lemm:goingthroughgeneralpoints} implies that the evaluation map $ev: \mathcal U \to X \times \mathbb P^1$ is dominant and separable. Moreover since it parametrizes sections, we have $\underline{\mathcal U} = M \times \mathbb P^1$. Then the composition $\mathcal U \to X \times \mathbb P^1 \to X$ with the projection is also dominant and separable. In particular, it is smooth at a general unmarked point $(x,y) \in \underline{\mathcal U} = M\times \mathbb{P}^1$ so that we have a surjection of tangents
 \[
 T_{M,x} \oplus T_{\mathbb{P}^1,y} \longrightarrow T_{X, ev(x,y)}
 \]
 whose kernel does not contain $ T_{\mathbb{P}^1,y}$.
 We may choose a general complete intersection $W \subset M$ through $x$ of dimension $\dim X - 1$.
 Then 
  \[
T_{W, x}\oplus T_{\mathbb{P}^1,y} \subset  T_{M,x} \oplus T_{\mathbb{P}^1,y} \longrightarrow T_{X, ev(x,y)}
 \]
 is still surjective so that the composition $\mathcal{C}_W \to X\times \mathbb{P}^1 \to X$ is separable and dominant. Thus our assertion follows.

(2) Suppose that our log variety $(X, \Delta)$ is separably $\varsigma$-rationally connected. This will imply that 
there is a family $\pi : \mathcal U \to W$ of log rational curves with contact orders $\varsigma$ and evaluation map $s : \mathcal U \to X$ such that $\underline{\mathcal U} = W \times \mathbb P^1$ and the evaluation map $ev^{(2)} : \mathcal U\times_W \mathcal U \to X^2$ is dominant and separable.  We use a similar construction in (1) but instead of using $\mathbb G_m$ we use the $2$-dimensional affine group $\mathbb G_a \rtimes \mathbb G_m$. Thus we construct $\widetilde{\mathcal U} \to \widetilde{W} = (\mathbb{G}_{a}\rtimes \mathbb G_m) \times W$ 
with the evaluation map $\widetilde{s} : \widetilde{\mathcal U} \to X \times \mathbb P^1$ using multiplication by a $\mathbb{G}_{a}\rtimes \mathbb G_m$-action on $\mathbb{P}^{1}$.
Then $\widetilde{ev}^{(2)} : \widetilde{\mathcal U}\times_{\widetilde{W}} \widetilde{\mathcal U} \to (X \times \mathbb P^1)^2$ is dominant and separable. 
This will imply that there is a component $M$ of $\Sec_{\varsigma}(X\times \mathbb P^1/\mathbb P^1)$ with universal family $\pi : \mathcal U \to M$ such that $ev^{(2)} : \mathcal U\times_M \mathcal U \to (X\times \mathbb P^1)^2$ is dominant and separable. Then it follows from Lemma \ref{lemm:goingthroughgeneralpoints} that a general $f : C \to X \times \mathbb P^1$ parametrized by $M$ is very free.

Conversely suppose that there is a very free log section $f : C \to X \times \mathbb P^1$ with contact orders $\varsigma$. Let $M$ be a component of $\Sec_{\varsigma}(X\times \mathbb P^1/\mathbb P^1)$ containing $f$. Let $\pi : \mathcal U \to M$ be the universal family over $M$ with the evaluation map $ev : \mathcal U \to X \times \mathbb P^1$. Then since $N_f$ is ample, for any points $p, q$ on $C$, $H^0(C, N_f) \to N_f|_p \oplus N_f|_q$ is surjective.  This means that for general $p, q$, we have a surjection $H^0(C, N_f) \oplus T_C|_p \oplus T_C|_q \to f^*T_{X \times \mathbb P^1}|_p \oplus f^*T_{X \times \mathbb P^1}|_q$. We conclude that the evaluation map $ev^{(2)} : \mathcal U \times_M \mathcal U \to (X \times \mathbb P^1)^2$ is dominant and separable proving the claim.
\end{proof}

\begin{proof}[Proof of Proposition \ref{prop:uniruled=free}]
This follows directly from Lemmas \ref{lem:free-curve=free-section} and \ref{lemm:varsigma-rat_connected}. 
\end{proof}

\section{Campana maps and curves}
\label{sec:campanamaps}

Here we recall the definitions of Campana curves and Campana sections. First let us define Campana orbifolds and Campana fibrations:

\begin{defi}
Let $\underline{X}$ be a smooth projective variety and $\Delta = \sum_i \Delta_i$ be a strict normal crossings divisor on $\underline{X}$ such that $\Delta_i$ is irreducible. Let $X$ be the log scheme associated to the pair $(\underline{X},\Delta)$. 
For each $i$ let $m_i$ be either a positive integer $\geq 1$ or $\infty$ and set $\epsilon_i = 1- \frac{1}{m_i}$. We define 
\[
\Delta_\epsilon = \sum_i \left(1- \frac{1}{m_i}\right)\Delta_i
\]
and call the pair $(X, \Delta_\epsilon)$ a {\em Campana orbifold}. We say $(X, \Delta_\epsilon)$ is a {\em klt Campana pair} if all $m_i$ are positive integers.  This is equivalent to saying that the pair $(\underline{X}, \Delta_\epsilon)$ has only klt singularities.

A Campana fibration over $B$ is a Campana orbifold $({\mathcal X}, \Delta_\epsilon)$ with a log fibration $\pi : (\mathcal X, \Delta) \to B$. Similarly we say $(\mathcal X, \Delta_\epsilon)/B$ is a {\em klt Campana fibration} if all the $m_i$ are positive integers.

\end{defi}

Next we define the notion of Campana maps:

\begin{defi}
\label{defi:Campanasection}
Let $(X, \Delta_\epsilon)$ be a Campana orbifold.  Suppose $\varsigma$ is a collection of positive contact orders as in \eqref{eq:discrete-data} for stable log maps to $X$.  As before we will let $c_{k,i}$ denote the multiplicity of the $i$th irreducible component $\Delta_{i}$ along the $k$th marked point.

We say that $\varsigma$ is of {\em Campana type} if for every irreducible component $\Delta_{i}$ of $\Delta$: 
\begin{enumerate}
 \item when $m_i < \infty$, every index $k$ satisfies either $c_{k,i} = 0$ or $c_{k, i} \geq m_i$. 
 \item when $m_i = \infty$, there is at most one index $k$ with $c_{k,i} > 0$. 
\end{enumerate}
A stable log map $f \colon C \to X$ over a log point $S$ is called {\em a Campana map} if the collection of its contact orders are of Campana type. 
We call $f$ a {\em Campana curve} if $f$ is non-degenerate, or equivalently $\mathcal{M}_S$ is trivial.
Recall that if $f$ is non-degenerate, then $C$ is smooth and $f^{-1}(\Delta)$ consists of only marked points.

Next let $\pi \colon (\mathcal X, \Delta_\epsilon) \to B$ be a Campana fibration over a smooth projective curve $B$ of genus $g$.
A genus $g$ stable log map $f \colon C \to \mathcal X$ over a geometric log point $S$ is called a {\em stable Campana map of section type} if $f$ is a Campana map to $(\cX,\Delta_\epsilon)$, and the composition $C \to \cX \to B$ has degree $1$.
When $C$ is a Campana curve, we call $f$ a {\em Campana section}.
\end{defi}

The condition $(2)$ means that when $m_i = \infty$ and $f(p_{k}) \in \Delta_i$, then $f(C\setminus \{p_k\}) \cap \Delta_i = \emptyset$. In particular, when every multiplicity is $\infty$ then the $\mathbb{A}^1$-curves in $X$ studied in \cite{CZ17} are examples of Campana curves.

\begin{rema}
In Definition~\ref{defi:Campanasection}, the condition on contact orders when $m_i = \infty$ is slightly different than the typical definition used by arithmetic geometers.  The exactly analogous definition is to fix a finite set $S$ of places on $B$ and to insist that the intersections of $\Delta_i$ and $f(C)$  can only happen above these finitely many points.
Our definition allows for more flexibility.
\end{rema}

\begin{rema}
An alternative method for dealing with infinite contact orders in Definition~\ref{defi:Campanasection} would be to insist that our curve meets $\cup_{i | m_{i} = \infty} \Delta_{i}$ at only one point.  This would be more closely analogous to the notion of an $\mathbb{A}^{1}$-curve.  It has the additional advantage that if we construct $\Delta'_{\epsilon}$ from $\Delta_{\epsilon}$ by reducing the multiplicities then a Campana curve for $(X,\Delta_{\epsilon})$ can be transformed into a Campana curve for $(X,\Delta'_{\epsilon})$ by taking a finite cover ramified at the point of intersection with the multiplicity $\infty$ divisors.  On the other hand, this alternative definition is not as flexible as Definition~\ref{defi:Campanasection}.
\end{rema}

\begin{rema} \label{rema:birationaldependence}
Suppose that $(X,\Delta_\epsilon)$ is a Campana orbifold and that $\phi: X' \to X$ is a birational morphism.  It is natural to equip $X'$ with a Campana orbifold structure $(X',\Delta'_\epsilon)$ in a ``minimal'' way (see \cite[Section 3.6]{PSTVA}).  However, with this choice the strict transform of a Campana curve on $X$ need not be a Campana curve on $X'$.  In fact, there does not seem to be a natural way of defining a ``birationally invariant'' theory of Campana curves.

Next suppose that $\pi: (\mathcal{X},\Delta_\epsilon) \to B$ is a Campana fibration.  The previous paragraph shows that the set of Campana sections can really depend on the integral model $\mathcal{X}$ and not just on the generic fiber $\mathcal{X}_{K(B)}$.  For this reason we will always specify the integral model $(\mathcal{X},\Delta)$ when discussing Campana sections.
\end{rema}

\subsection{Campana uniruledness and Campana rational connectedness}

Here we introduce the notation of (separable) Campana uniruledness and (separable) Campana rational connectedness:

\begin{defi}
Let $(X, \Delta_\epsilon)$ be a Campana orbifold.
\begin{enumerate}
\item We say $(X, \Delta_\epsilon)$ is (separably) Campana uniruled if there exist a collection of contact orders $\varsigma$ of  Campana type such that $X$ is (separably) $\varsigma$-uniruled. Furthermore, when all non-zero contact orders of $\varsigma$ are not divisible by $\mathrm{char} \, \mathbf k$, we say $(X, \Delta_\epsilon)$ is (separably) Campana uniruled by good contact orders.
\item We say $(X, \Delta_\epsilon)$ is (separably) Campana rationally connected if there exist a collection of contact orders $\varsigma$ of Campana type such that $X$ is (separably) $\varsigma$-rationally connected. We also define (separably) Campana rational connectedness by good contact orders in an analogous way.
\end{enumerate}
\end{defi}

\begin{lemm}
Let $(X, \Delta)$ be a Campana orbifold. 
\begin{enumerate}
\item A Campana orbifold $(X, \Delta_\epsilon)$ is separably Campana uniruled by good contact orders if and only if there is a non-trivial free Campana section of contact orders $\varsigma$ on $(X\times \mathbb P^1/\mathbb P^1, \Delta_\epsilon\times \mathbb P^1)$ such that every non-zero contact order of $\varsigma$ is not divisible by $\mathrm{char}\, \mathbf k$; 
\item a Campana orbifold $(X, \Delta_\epsilon)$ is separably Campana rationally connected by good contact orders if and only if there is a very free Campana section of contact orders $\varsigma$ on $(X \times \mathbb P^1/\mathbb P^1, \Delta_\epsilon\times \mathbb P^1)$ such that every non-zero contact order of $\varsigma$ is not divisible by $\mathrm{char}\, \mathbf k$.
\end{enumerate}
\end{lemm}
\begin{proof}
This follows from Lemma~\ref{lemm:varsigma-rat_connected}. 
\end{proof}

An important conjecture in this direction, due to Campana, is: 

\begin{conj}[Campana]
Assume that our ground field has characteristic $0$.
Let $(X, \Delta_\epsilon)$ be a klt Fano orbifold, i.e., $(X, \Delta_\epsilon)$ is a klt Campana orbifold such that $-(K_{\underline{X}} + \Delta_\epsilon)$ is ample. Then $(X, \Delta_\epsilon)$ is Campana rationally connected. 
\end{conj}

For our applications, we will need the following conjecture.

\begin{defi}
\label{defi:stronglycampanauniruled}
Let $(X, \Delta_\epsilon)$ be a klt Campana orbifold. We say that $(X, \Delta_\epsilon)$ is strongly Campana uniruled if there exists a free Campana curve $f : C \to X$ such that the class $f_{*}[C]$ is in the interior of the nef cone $\mathrm{Nef}_1(X)$ of curves. 
\end{defi}

\begin{conj}
\label{conj:main}
Assume that our ground field has characteristic $0$.
Then every klt Fano orbifold $(X, \Delta_\epsilon)$ is strongly Campana uniruled.
\end{conj}

Note that when $X$ has Picard rank $1$, being in the interior of the nef cone is automatic.

\begin{prop}
\label{prop:codimension1strata}
Let $(X, \Delta_\epsilon)$ be a klt Fano orbifold. Assume that Conjecture~\ref{conj:main} holds. Then for each $\Delta_i$, there exists a free Campana curve $f : C \to X$ such that $f(C)$ meets with the codimension $1$ stratum of $\Delta_i$. 
\end{prop}

\begin{proof}
This follows from Proposition~\ref{prop:splittingcontactorder} and Conjecture~\ref{conj:main}.
Indeed, one can split contact orders inductively as described in Section~\ref{ss:contact-splitting}.
While doing so, freeness will be preserved by Remark~\ref{rema:splittingcontact_freeness}.
\end{proof}

\section{Weak approximation}
\label{sec:weakapproximation}

Weak approximation for Campana sections looks somewhat different than weak approximation for sections.  Indeed, by Remark \ref{rema:birationaldependence} the notion of a Campana section depends on the choice of integral model $\pi: (\mathcal{X},\Delta) \to B$ and not just the generic fiber.  

\subsection{Setting up weak approximation for Campana pairs}

Let $B$ be a smooth projective curve over an algebraically closed field.  Suppose that $(\mathcal{X}_{\eta},\Delta_{\eta})$ is a smooth projective geometrically integral klt Campana pair over the function field $K$ of $B$.  For any place $b$ of $K$, we denote by $\widehat{\mathcal{O}}_{B, b}$ the completion of the local ring $\mathcal{O}_{B, b}$ with respect to the maximal ideal and by $K_{b}$ the fraction field of $\widehat{\mathcal{O}}_{B, b}$.

\begin{defi}
Let $\pi: (\mathcal{X},\Delta_\epsilon) \to B$ be a klt Campana fibration.
We say that $b \in B$ is a place of good reduction if 
there exists a regular model $\pi' : \underline{\mathcal X}' \to \Spec \, \mathcal O_{B, b}$ of the generic fiber $\underline{\mathcal X}_\eta$ such that the special fiber $\underline{\mathcal X}'_b$ is smooth.
\end{defi}

Note that this definition is the same as the usual one and does not rely on the log structure in any way.

\begin{defi}
\label{defi:Campanajets}

Let $\pi: (\mathcal{X},\Delta_\epsilon) \to B$ be a klt Campana fibration.
Let $\Spec \, \mathbf k[t]/(t^{n+1})$ be the $n$-th jet scheme.  We say an $n$-th jet $\sigma : \Spec \, \mathbf k[t]/(t^{n+1}) \to \underline{\mathcal X}$ is an {\em admissible} $n$-th jet if the composition $$\Spec \, \mathbf k[t]/(t^{n+1}) \to \underline{\mathcal X} \to B$$ is a closed embedding.

Let $I_{\sigma}$ denote the set of indices $i$ such that $\sigma(\Spec \, \mathbf k) \subset \Delta_{i}$.  We say an admissible $n$-th jet $\sigma$ is a {\em Campana} $n$-th jet if it satisfies:
\begin{enumerate}
\item $n \geq \max_{i} \{m_{i}\}_{i \in I_{\sigma}}$, and
\item for $i \in I_{\sigma}$ the ideal defined by the pullback of $\Delta_i$ is given by $(t^m)$ with $m \geq m_i$.
\end{enumerate}

We say a Campana fibration  $\pi: (\mathcal{X},\Delta_\epsilon) \to B$ satisfies weak approximation if for any finite number of Campana jets sitting in distinct fibers, there is a Campana section $f : C \to \mathcal X$ which induces the given Campana jets.  Similarly, we say that $\pi$ satisfies weak approximation at places of good reduction if any finite number of Campana jets sitting in distinct fibers of good reduction are induced by a Campana section. 
\end{defi}

\subsection{Deformation theory of log sections while fixing jets}
In Section \ref{ss:point-constraints-deformation} we discussed the deformation theory of log sections through fixed points.  In this section we extend these results to discuss sections through fixed jets.  The strategy is the usual one (e.g., \cite[Section 2.3]{HT06}): we blow-up to translate jet data into incidence data.

Suppose that $\pi: (\mathcal X, \Delta_\epsilon) \to B$ is a klt Campana fibration. Let $\{(p_j, \sigma_j)\}$ be a finite set of admissible jets living in distinct fibers. 
We assume that these jets are not supported on the boundary $\Delta$.
 Suppose that the jet $\sigma_{j}$ lies over the point $b_{j} \in B$.  Let $\widehat{B}_{b}$ denote the completion of $B$ at $b$ and let $\widehat{\pi}: \widehat{\mathcal{X}}_{b} \to \widehat{B}_{b}$ denote the base-change of $\pi$.  By Hensel's Lemma, each $\sigma_{j}$ is induced by a jet $\widehat{\sigma}_{j}$ over $\widehat{B}_{b_{j}}$. 

We then replace $\mathcal{X}$ with the following birational modification.  For each jet $\sigma_{j}$, we repeatedly perform point blow-ups in the fiber over $b_{j}$ where at each step we blow-up the point defined by the strict transform of $\widehat{\sigma}_{j}$.  The result will be a smooth birational model $\mathcal{X}'$ of $\mathcal{X}$ equipped with a morphism $\pi': \mathcal{X}' \to B$ with the following property: for each $j$, there is a distinguished irreducible component $E_{j}$ of the fiber over $b_{j}$ such that a section $C$ of $\pi$ will be tangent to $\sigma_{j}$ if and only if the strict transform of $C$ meets $E_{j}$.  We say that $\phi: \mathcal{X}' \to \mathcal{X}$ extracts the jets $\{ \sigma_{j} \}$.
For an $n$-th admissible jet, one can construct such a $\mathcal X'$ by blowing up $n+1$ times.

\begin{lemm}
Let $\pi: (\mathcal{X}, \Delta_\epsilon) \to B$ be a klt Campana fibration.  Suppose $\{ \sigma_{j} \}_{j=1}^{r}$ is a finite set of admissible jets supported on distinct fibers such that the support of jets does not lie on $\Delta$.  
Let $\phi: \mathcal{X}' \to \mathcal{X}$ be the birational model extracting the jets $\{ \sigma_{j} \}$.  Suppose that $C$ is a log section of $\pi$ with the canonical log structure that approximates the jets $\{ \sigma_{j} \}$ and let $C'$ be its strict transform on $\mathcal{X}'$.  Then
\begin{equation*}
N_{C'/\mathcal{X}'} \cong N_{C/\mathcal X}\left(-\sum n_jp_j\right)
\end{equation*}
In particular the deformation theory of sections $C$ of $\pi$ approximating $\{ \sigma_{j} \}_{j=1}^{r}$ is given by the twisted normal bundle in the equation above.

\end{lemm}

\begin{proof}
Recall that $\phi$ is a composition of blow-ups at smooth points.  Arguing inductively, it suffices to compute what happens when $\phi$ is the blow-up of a single point $p$ not contained in $\Delta$.  Letting $i: E \to \mathcal{X}'$ denote the inclusion of the exceptional divisor with the trivial log structure, we have an exact diagram of sheaves
\begin{equation*}
\xymatrix{
& 0 \ar[d] & 0 \ar[d] & 0 \ar[d] & \\
0 \ar[r] & T_{{C}'} \ar[r] \ar[d]^{=} & T_{{\mathcal{X}}'}|_{{C}'} \ar[r] \ar[d] & N_{{C}'/{\mathcal{X}}'} \ar[r] \ar[d] & 0 \\
0 \ar[r] & T_{{C}'} \ar[r] \ar[d] & \phi^{*}T_{{\mathcal{X}}}|_{{C}'} \ar[r] \ar[d] & N_{{C}/{\mathcal{X}}} \ar[r] \ar[d] & 0 \\
& 0 \ar[r]  & i_{*}T_{{E}}(E)|_{{C}'} \ar[r]_(.6){=} \ar[d] & \mathcal{K} \ar[r] \ar[d] & 0 \\
& & 0 & 0 & 
},
\end{equation*}
where $\mathcal K$ is the cokernel of $N_{{C}'/{\mathcal{X}}'} \to N_{{C}/{\mathcal{X}}}$.
We conclude that $N_{{C}'/{\mathcal{X}}'} \cong N_{{C}/{\mathcal X}}(-p)$ as desired.
\end{proof}

\subsection{Weak approximation in characteristic $0$}

Our goal in this section is to prove the following theorem:

\begin{theo}
\label{theo:weakapproximation}
Assume that $\mathbf k$ is an algebraically closed field of characteristic $0$.
Let $S$ be a finite set of closed points on $B$ and we assume that $\underline{\pi} : \underline{\mathcal X} \to B$ satisfies weak approximation outside of $S$.
Let $\pi: (\mathcal{X},\Delta_\epsilon) \to B$ be a klt Campana fibration over $\mathbf k$ such that a general fiber of $\pi$ is rationally connected and is strongly Campana uniruled. Suppose we fix points $p_{1},\ldots,p_{r} \in \mathcal{X}$ supported on distinct fibers outside of $S$
and for each index $j$ we choose a Campana $n_j$-th jet $\sigma_{j}$ at $p_{j}$. 
Then there exists a Campana section approximating the jet data $\{ p_{j}, \sigma_{j} \}$.
\end{theo}

We assume that $p_1, \cdots, p_s$ are contained in $\Delta$ and $p_{s+1}, \cdots, p_r$ are not contained in $\Delta$.
First we prepare the following lemma:

\begin{lemm}
\label{lemm:weakapproximationbyfreesection}
Assume that $\mathbf k$ is an algebraically closed field of characteristic $0$.
Let $S$ be a finite set of closed points on $B$ and we assume that $\underline{\pi} : \underline{\mathcal X} \to B$ satisfies weak approximation outside of $S$.
Let $\underline{\pi} : \underline{\mathcal X} \to B$ be a flat projective morphism from a smooth projective variety such that a general fiber is rationally connected and let $\rho \geq 0$ be any non-negative real number.
For $j = 1, \cdots, r$, let $\{p_j, \sigma_j\}$ be an admissible $n_j$-th jet such that $p_1, \cdots, p_r \in \underline{\mathcal X}$ are sitting in fibers outside of $S$.
Then $\{\sigma_j\}$ can be approximated by a section $f: \underline{C} \to \underline{\mathcal X}$ of $\underline{\pi}$ such that 
\[
\mu^{min}\left(N_{\underline{C}/\underline{\mathcal X}}\left(-\sum_j n_jp_j\right)\right)\geq \rho + 2g(B).
\]
\end{lemm}

\begin{proof}
First let us assume that $\mathbf k$ is uncountable.
We fix $\ell = \lceil \rho \rceil + 2g(B) + 1$ very general points $q_1, \cdots, q_\ell$. 
By our assumption, we can find a section $\underline{f}' : \underline{C}' \to \underline{\mathcal X}$ that approximates our finite number of jets $\{\sigma_j\}$ and also goes through our $\ell$ very general points. 
Let $M$ be the irreducible component of the moduli space parametrizing sections which approximate our jets $\{\sigma_j\}$ such that $\underline{f}' \in M$.  
Since our points are very general, a general section $\underline{C}$ parametrized by $M$ goes through $\ell$ very general points.
This implies that $N_{\underline{C}/\underline{\mathcal X}}(-\sum_j n_jp_j-\sum_{k = 1}^{\ell-1}q_k)$ is generically globally generated. Thus it follows from \cite[Lemma 2.6]{LRT23} that we have
\[
\mu^{min}\left(N_{\underline{C}/\underline{\mathcal X}}\left(-\sum_j n_jp_j\right)\right)\geq \lceil \rho \rceil + 2g(B)
\]
verifying our assertion.

When $\mathbf k$ is countable, we consider the moduli space $M$ of sections approximating the jet data  $(p_j, \sigma_j)$. Let $\mathbf k' \supset \mathbf k$ be an uncountable algebraically closed field. Then there exists an irreducible component $M_1 \subset M$ such that $M_1\otimes \mathbf k'$ contains an open subset $U$ parametrizing sections $\underline{C}$ such that
\[
\mu^{min}\left(N_{\underline{C}/\underline{\mathcal X}}\left(-\sum_j n_jp_j\right)\right)\geq \lceil \rho \rceil + 2g(B).
\]
Since $\mathbf k$-valued points on $M_1$ are Zariski dense, one can find a section $\underline{C}$ defined over $\mathbf k$ that satisfies 
\[
\mu^{min}\left(N_{\underline{C}/\underline{\mathcal X}}\left(-\sum_j n_jp_j\right)\right)\geq \lceil \rho \rceil + 2g(B).
\]
Thus our assertion follows.
\end{proof}

We also need another lemma:

\begin{lemm}
\label{lemm:birationalmodification}
Let $\underline{\pi} : \underline{\mathcal X} \to B$ be a flat projective morphism from a smooth projective variety with a boundary divisor $\Delta = \cup_i \Delta_i$ which is a SNC divisor and flat over $B$. For $j = 1, \cdots, r$, let $\{p_j, \sigma_j\}$ be an admissible $n_j$-th jet such that each point $p_1, \cdots, p_r \in \underline{\mathcal X}$ is a smooth point of the fiber containing it. Suppose that we have a section $C$ approximating these jets and not contained in the support of $\Delta$. Then there is a birational morphism $\beta: \underline{\widetilde{\mathcal X}} \to \underline{\mathcal X}$ from a smooth projective variety $\underline{\widetilde{\mathcal X}}$ such that the strict transform $\widetilde{\Delta}$ of $\Delta$ is SNC, $\beta$ is an isomorphism over $\underline{\mathcal X} \setminus \{p_1, \cdots, p_r\}$, and the strict transform $\widetilde{C}$ of $C$ does not meet with $\widetilde{\Delta}$ over $\underline{\pi}(p_j)$ for any $j = 1, \cdots, r$.
\end{lemm}

\begin{proof}
It is enough to prove the statement for one section $C$ and a point $p \in C$. If $p \not\in \Delta$, then there is nothing to prove. Suppose that $p \in \Delta$. Then we successively blow up the support $p$ of the strict transform of a section $f : C \to \underline{\mathcal X}$.  
Then the local intersection multiplicity of the strict transform $\widetilde{C}$ and the strict transform $\widetilde{\Delta}$ is strictly decreasing along successive blow-ups, so eventually it becomes $0$, proving the claim. Note that the smoothness of $p$ in its fiber is also preserved due to the fact that it is the intersection of the fiber with a section. Also note that if the local intersection of $C$ and $\Delta_i$ is given by $k_i$, then we need to perform successive blow ups $\max \{ k_i\}$-times. 
\end{proof}

Here is our strategy for proving Theorem \ref{theo:weakapproximation}.  Suppose we fix a finite set of Campana jets supported on fibers of good reduction.  By the assumption we can find a section $\underline{C}$ which induces this finite set of jets.  Of course $\underline{C}$ might not satisfy the Campana condition at the other points of intersection $\{q_{k}\}$ with $\Delta$.  By gluing on $\pi$-vertical Campana curves and smoothing (while leaving the jets fixed), we increase the intersection numbers at the $q_{k}$ to achieve the Campana condition everywhere.

\begin{proof}[Proof of Theorem~\ref{theo:weakapproximation}]

Let $\underline{C}$ be a general section in $M$ obtained in the proof of Lemma~\ref{lemm:weakapproximationbyfreesection} with $\rho = 2g$.
If at a place $b = \pi(p)$ of a Campana $n$-th jet $(p, \sigma)$, the local intersection multiplicity of $C$ and $\Delta_i$ is greater than $n$ for some $i$, then we may replace $(p, \sigma)$ by a deeper jet of $C$ at $p$ so that the local intersection multiplicity of $C$ and $\Delta_i$ is smaller than $n$ for any $i$.
Let $\beta : \underline{\widetilde{\mathcal X}} \to \underline{\mathcal X}$ be a birational projective morphism constructed in Lemma~\ref{lemm:birationalmodification}.
By taking the strict transform, we have a stable map $\underline{f} : \underline{C} \to \underline{\widetilde{\mathcal X}}$. Then note that $\underline{C}$ meets with the exceptional divisor of $\beta$ only over points $\pi(p_j)$.
We impose the log structure on $\underline{\widetilde{\mathcal X}}$ associated to the pair $(\underline{\widetilde{\mathcal X}}, \widetilde{\Delta})$ where $\widetilde{\Delta}$ is the strict transform of $\Delta$. We can think of $(\widetilde{\mathcal X},\widetilde{\Delta}_{\epsilon})$ as a Campana pair by equipping each irreducible component of $\widetilde{\Delta}$ with the same multiplicity as the corresponding component of $\Delta$.
Then we impose the minimal log structure on $C \to S$ so that $f : C \to \widetilde{\mathcal X}$ is a log section with only contact markings. Let $\{ q_k \}$ denote the set of marked points on $C$. By construction $f(q_k)$ is a smooth point of $\widetilde{\Delta}$ and $f(C)$ meets with $\widetilde{\Delta}$ transversally at those points. Let $\ell$ be the number of $q_k$'s where the Campana condition is not satisfied.

It follows from the proof of Lemma~\ref{lemm:generalpoints} that $f : C \to \widetilde{\mathcal X}$ goes through $2g+1$ general points while approximating the jet data $\{ (\widetilde{p}_{j}, \widetilde{\sigma}_{j}) \}$ induced by $\{(p_j, \sigma_j)\}$.
In particular, this implies that
\[
\mu^{min}\left(N_{f}\left(-\sum_{j = 1}^r \widetilde{n}_j \widetilde{p}_j\right)\right) \geq 2g.
\]
This means that $N_{f}(-\sum_{j = 1}^r \widetilde{n}_j \widetilde{p}_j)$ has vanishing $H^1$ and is globally generated. Let $q_k$ be a marked point on $C$ which does not satisfy the Campana condition for $(\widetilde{\mathcal X},\widetilde{\Delta}_{\epsilon})$.
By the construction, $f(q_k)$ is a general point in a codimension $1$ strata of a general fiber.  By assumption this fiber possesses a free Campana rational curve which is in the interior of the nef cone.
Proposition~\ref{prop:codimension1strata} implies that we have a free $\pi$-vertical Campana rational curve passing through $f(q_k)$. 
Then we glue $C$ and $T$ via a contracted component $L$ so that we obtain a glued stable log map $\widetilde{f} : \widetilde{C} \to \widetilde{\mathcal X}$ using the construction of Section~\ref{subsec:gluing}.
We claim that Proposition~\ref{prop:smoothing lemma} allows us to smooth $\widetilde{f} : \widetilde{C} \to \widetilde{\mathcal X}$ while approximating the jet data $\{ \widetilde{p}_{j}, \widetilde{\sigma}_{j} \}$. 
Indeed, the assumptions of Proposition~\ref{prop:smoothing lemma} follow from Lemma~\ref{lemm:N-at-marking} and Proposition~\ref{prop:deformation-relative-boundary}.
We denote its general smoothing as a log section $f_1 : C_1 \to \widetilde{\mathcal X}$. Then the number of marked points $q_k$ which does not satisfy Campana condition is $\ell-1$. Moreover since $f : C \to \widetilde{\mathcal X}$ goes through $2g+1$ general points, $f_1 : C_1 \to \widetilde{\mathcal X}$ goes through $2g + 1$ general points as well. Thus by repeatedly applying the construction of Section~\ref{subsec:gluing} to marked points which do not satisfy the Campana condition and smoothing the resulting stable maps we can construct a Campana log section $f_\ell : C_\ell \to \widetilde{\mathcal X}$ which approximates the jet data $\{ \widetilde{p}_{j}, \widetilde{\sigma}_{j} \}$. Then $h: C_\ell \to \underline{\widetilde{\mathcal X}} \to \underline{\mathcal X}$ witnesses our assertion.
\end{proof}

%\begin{rema}
%\label{rema:weakapproximation}
%The proof of Theorem~\ref{theo:weakapproximation} shows more. When the generic fiber $\underline{\mathcal X}_\eta$ satisfies the usual weak approximation, our proof actually shows that one can find a Campana section approximating Campana jets at any finite set of places including places of bad reduction.
%\end{rema}

\begin{proof}[Proof of Theorem~\ref{theo:mainI_intro}]
This follows from Theorem~\ref{theo:weakapproximation}.
\end{proof}

\begin{coro}
\label{coro:CRC}
Assume that $\mathbf k$ is an algebraically closed field of characteristic $0$.
Let $(X, \Delta_\epsilon)$ be a klt Campana orbifold such that $\underline{X}$ is rationally connected and $(X, \Delta_\epsilon)$ is strongly Campana uniruled. Then $(X, \Delta_\epsilon)$ is Campana rationally connected.
\end{coro}

\begin{proof}
First let us assume that $\mathbf k$ is uncountable. Let $\mathcal X = X \times \mathbb P^1$ and $\widetilde{\Delta}_\epsilon = \Delta_\epsilon \times \mathbb P^1$. We pick two very general points on $\mathcal X$. It follows from Theorem~\ref{theo:weakapproximation} that we can find a Campana section $C$ of $((\mathcal X, \widetilde{\Delta}_\epsilon)/\mathbb P^1)$ passing through these two points. Since these two points are very general, $C$ log-deforms to go through two general points on $\mathcal X$.  It follows from the proof of Lemma~\ref{lemm:generalpoints} that a general log-deformation of $C$ is very free. Thus our assertion follows. When $\mathbf k$ is countable, one can reduce to the statement for uncountable fields as in Lemma~\ref{lemm:weakapproximationbyfreesection}.
\end{proof}

\begin{proof}[Proof of Corollary~\ref{coro:CRC_intro}]
This follows from Corollary~\ref{coro:CRC}.
\end{proof}

\section{Campana curves for $\mathbb{P}^{1}$-fibrations}
\label{sec:P1}

Suppose that $(\mathbb{P}_{\eta}^{1},D_{\eta})$ is a klt Fano orbifold with a smooth integral model $\pi: S \to B$.  Explictly, this means that:
\begin{itemize}
\item $S$ is smooth and comes equipped with a SNC divisor $\mathcal{D}$ that is flat and generically SNC over $B$.
\item The general fiber of $\pi$ is isomorphic to $\mathbb{P}^{1}$.
\end{itemize}
The possible coefficient choices depend on the support of $\mathcal{D}$:
\begin{enumerate}
\item If $\mathcal{D}$ consists of a single section or degree $2$ multisection, we can assign to it any integer $m \geq 2$.
\item If $\mathcal{D}$ consists of two sections, we can assign to the pair any integers $m_{1},m_{2} \geq 2$.
\item If $\mathcal{D}$ consists of a single degree $3$ multisection, we must assign it the multiplicity $m = 2$.
\item If $\mathcal{D}$ consists of a section $D_1$ and a degree $2$ multisection $D_2$, we must assign them the multiplicities $(m_1, m_2) = (m, 2)$ or $(2, 3)$ where $m \geq 2$ is any integer.
\item If $D$ consists of three sections, then multiplicities are $(m_1, m_2, m_3) = (2, 2, m), (2, 3, 3), (2, 3, 4)$ or $(2, 3, 5)$ where $m$ is any integer $m \geq 2$.
\end{enumerate}

Our goal of this section is to show Conjecture~\ref{conj:main} for $(\mathbb{P}^{1},D)$ where $D = \mathcal{D}|_{\mathbb{P}^{1}}$ for a general fiber of $\pi$.

\subsection{The existence of free Campana rational curves in the absolute case}
\label{subsec:step1}

Assume that $(\mathbb P^1, \Delta_\epsilon)$ is a klt Fano orbifold. When we have two Campana orbifolds $(\mathbb P^1, \Delta_\epsilon)$ and $(\mathbb P^1, \Delta'_\epsilon)$ such that $\Delta_\epsilon \geq \Delta'_\epsilon$, then a Campana curve with respect to $\Delta_\epsilon$ is automatically a Campana curve with respect to $\Delta'_{\epsilon}$. So we may assume that 
\begin{enumerate}
\item $\Delta$ consists of two points with $(m_1, m_2) = (m, m)$ where $m \geq 2$ is an integer;
\item $\Delta$ consists of three points with $(m_1, m_2, m_3) = (2, 2, m)$ where $m \geq 2$ is an integer, or;
\item $\Delta$ consists of three points with $(m_1, m_2, m_3) = (2, 3, 5)$.
\end{enumerate}

Using this simplification we show:
\begin{theo} \label{theo:orbifoldp1}
Let $(\mathbb P^1, \Delta_\epsilon)$ be a klt Fano orbifold. Then there exists a stable log map $f : C \to \mathbb P^1$ of genus $0$ such that $f$ is a free Campana curve. Moreover, one can choose $f$ so that the log normal sheaf of $f$ has degree $m$ where $m$ is any non-negative integer.  

\end{theo}

\begin{proof}
We first assume that our ground field $\mathbf{k}$ has characteristic $0$.  In the case (1), one can find a degree $m$ cover $f ': \mathbb P^1 \to \mathbb P^1$ totally ramified at the support of $\Delta$.  If we equip the domain with the log structure associated to $f'^{-1}(\Delta)$ then $f'$ is a Campana rational curve.  It is also free since the log normal sheaf has degree $0$.  
Let $g : \mathbb P^1 \to \mathbb P^1$ be a general degree $d$ cover and equip the domain with the log structure associated to $g^{-1}f'^{-1}(\Delta)$. Then $f = f'\circ g$ is a Campana curve and the log normal sheaf has degree $2d -2$.  This achieves every non-negative even degree for the log normal sheaf; for an odd degree, we can instead let $g$ be simple ramified at one point of $g^{-1}f'^{-1}(\Delta)$. Thus our assertion follows. 

In the cases (2) and (3), we claim that there exists a cover $f' : \mathbb{P}^{1} \to \mathbb P^1$ branched at the support of $\Delta$ with branch datum $$(m_1, \cdots, m_1), (m_2, \cdots, m_2), (m_3, \cdots, m_3).$$  In fact, these are classical examples of Belyi maps obtained by taking quotients of $\mathbb{P}^{1}$ by finite subgroups of $\mathrm{PGL}_{2}(\mathbf{k})$.  Case (2) corresponds to quotients by a dihedral group of order $2m$ whose action has stabilizers of size $2,2,m$.  This gives a degree $2m$ map $f': \mathbb{P}^{1} \to \mathbb{P}^{1}$ with branch data $(2^{m}), (2^{m}), (m,m)$.  Case (3) corresponds to a quotient by the icosahedral group whose action has stabilizers of size $2,3,5$.  This gives a degree $60$ map $f': \mathbb{P}^{1} \to \mathbb{P}^{1}$ with branch data $(2^{30}), (3^{20}), (5^{12})$.
It is clear that in each case $f'$ defines a Campana curve and the degree of log normal sheaf is given by $0$. Let $g : \mathbb P^1 \to \mathbb P^1$ be a general degree $d$ cover. Then $f = f'\circ g$ is a Campana curve and the log normal sheaf has degree $2d -2$. 
For odd degrees, we may let $g$ simply ramified at one of $g^{-1}f'^{^1}(\Delta)$. Thus our assertion follows.

We next assume that our ground field $\mathbf{k}$ has characteristic $p$.  In the case (1), we can use the same argument after perhaps increasing the degree of the cover to ensure that it is coprime to $m$.  In the cases (2) and (3), we can again look for a cover $f': \mathbb{P}^{1} \to \mathbb{P}^{1}$ constructed by taking a quotient by a finite subgroup of $\mathrm{PGL}_{2}(\mathbf{k})$.  The classification of such subgroups is presented in \cite[Theorem B, Theorem C, Remark 2.3]{Faber23}: in every characteristic $\mathrm{PGL}_{2}(\mathbf{k})$ admits subgroups isomorphic to $A_{5}$ and to dihedral subgroups of all orders.  For $A_{5}$, suppose the non-trivial stabilizers for the group action have size $s_{1},\ldots,s_{k}$.  Since this group scheme is reduced, we can apply Riemann-Hurwitz to the morphism $f': \mathbb{P}^{1} \to \mathbb{P}^{1}$ obtained by the quotient to see that
\begin{equation*}
-2 = 60 \cdot (-2) + \sum_{i=1}^{k} (s_{i}-1) \frac{60}{s_{i}}.
\end{equation*}
A quick numerical computation shows that we must have $\{ s_{1},s_{2},s_{3} \} = \{ 2,3,5\}$ or $\{ 2,2, 30 \}$.  The latter is ruled out by the explicit description in \cite[Proposition 4.21]{Faber23} in characteristic $\neq 5$ and by \cite[Proposition 4.22]{Faber23} in characteristic $5$.   
Thus the cover defines a free Campana curve.  For the dihedral group, for simplicity we may increase the size so that it has the form $2m$ with $m$ a prime different from $p$.  Repeating the Riemann-Hurwitz computation above, the stabilizers must have size $\{ 2,2, m \}$ or size $\{2m,2m\}$.  It suffices to show that in characteristic $2$ we can ensure that there are not two points that are fixed by the entire group.  This follows from \cite[Remark 2.1]{Faber23}.  The rest of the argument is the same as in characteristic $0$.

\end{proof}

\begin{rema} 
Consider the Campana orbifold $(\mathbb{P}^{1},\Delta_{\epsilon})$ with multiplicities $(2,3,4)$ in characteristic $2$.  \cite[Theorem A, Remark 2.3]{Faber23} shows that $\mathrm{PGL}_{2}(\mathbf{k})$ does not admit any subgroup isomorphic to the octahedral group $S_{4}$.  While we can still construct a free curve $(\mathbb{P}^{1},\Delta_{\epsilon})$ using an icosahedral subgroup instead, we do not know whether there exists a free curve cover $f': \mathbb{P}^{1} \to \mathbb{P}^{1}$ with the ``minimal'' possible ramification behavior.
\end{rema}

\begin{rema} \label{rema:campanasdef}
We can now compare Campana's definition of orbifold uniruledness and orbifold rational connectedness to our notions.  Suppose that $(X,\Delta_{\epsilon})$ is a smooth klt Campana orbifold and that $f: \mathbb{P}^{1} \to X$ is a morphism whose image is not contained in $\Supp(\Delta_{\epsilon})$.  In \cite{Campana11} Campana says that $f$ is an orbifold rational curve if we can give $\mathbb{P}^{1}$ an orbifold structure $(\mathbb{P}^{1},\Gamma_{\epsilon})$ such that:
\begin{enumerate}
\item $K_{\mathbb{P}^{1}} + \Gamma_{\epsilon}$ is antiample, and
\item for every point $p \in \mathbb{P}^{1}$ and every irreducible component $\Delta_{i}$ we have $t_{i,p} \geq \frac{m_{i,p}}{n_{i,p}}$ whenever $t_{i, p}$ is positive where $t_{i,p}$ is the local multiplicity of $\Delta_{i}$ along $p$ and $m_{i,p},n_{i,p}$ are respectively the multiplicities of $\Delta_{i}$ in $\Delta_{\epsilon}$ and $p$ in $\Gamma_{\epsilon}$.
\end{enumerate}
He then defines orbifold uniruledness and orbifold rational connectedness using the existence of families of orbifold rational curves through one or two general points respectively.

A priori Campana's notion of an orbifold rational curve is more general than our notion of a Campana curve.  Correspondingly, Campana's notions of uniruledness and rational connectedness are a priori more general than ours.  However, by precomposing a map $f: (\mathbb{P}^{1},\Gamma_{\epsilon}) \to (X,\Delta)$ by a free Campana curve $g: \mathbb{P}^{1} \to (\mathbb{P}^{1},\Gamma_{\epsilon})$ as in Theorem \ref{theo:orbifoldp1} we see that composition $f \circ g: \mathbb{P}^{1} \to (X,\Delta)$ is a Campana curve.  In particular, any family of orbifold rational curves through one (or two) general points yields a Campana curve through one (or two) general points, showing that the two notions are equivalent.  
\end{rema}

\section{Rational Campana curves in toric varieties}
\label{sec:toric}

\subsection{Log curves in toric targets}

\subsubsection{The targets}
Let $N \cong \bZ^{d}$ be a lattice and $M = N^{\vee}$ be its dual lattice. Write $N_{\mathbb{R}} = N\otimes_{\bZ}\bR$ 
and $M_{\mathbb{R}} = M\otimes_{\bZ}\bR$. 
For a fan $\Sigma$ in $N$, 
define the log scheme $X_{\Sigma}$ whose underlying variety is the toric variety $\underline{X_{\Sigma}}$ associated to the fan $\Sigma$ and whose log structure $\cM_{X_{\Sigma}}$ is associated to the toric boundary $\Delta_{\Sigma} \subset \underline{X_{\Sigma}}$.  Note that we do not require $\underline{X_{\Sigma}}$ to be smooth in general. However, the log scheme $X_{\Sigma}$ is always log smooth by \cite[Prop. 3.4]{KKato89}.
Let $\Sigma^{(1)} \subset \Sigma$ be the collection of rays. For a ray $\rho \in \Sigma^{(1)}$, denote by $u_\rho$ its lattice generator and let $\Delta_{\rho} \subset X_{\Sigma}$ be the corresponding torus invariant divisor.  Thus, we have $\Delta_{\Sigma} = \sum_{\rho \in \Sigma^{(1)}}\Delta_{\rho}$.

\subsubsection{Contact orders}\label{sss:toric-contact-orders}
Consider $X_{\Sigma}$ as the target of log maps. 
The set of possible contact orders at a marked point is in bijection with the lattice points in the support of the fan $|\Sigma|$.  In particular, if $\Sigma$ is complete, then a contact order at a marked point can be identified by an element of $N$.  Below we will explain this bijection when $\Sigma$ is a smooth fan; the general case is described by \cite[\S 5.2]{ACMW17} and \cite[\S 2.3.8]{ACGS20}.

In case $\Sigma$ is a smooth fan a lattice point $c \in N$ specifies a contact order as defined in \S \ref{sec:deformation} as follows. 
Let $\sigma \in \Sigma$ be the minimal cone containing $c$, and let $\sigma^{(1)}$ be the set of rays of $\sigma$. 
The smoothness of the cone $\sigma$ implies a unique presentation $\mathbf{c} = \sum_{\rho \in \sigma^{(1)}} c_{\rho} u_\rho$ for a positive integer $c_{\rho}$ by the minimality of $\sigma$. 
Thus $\mathbf{c}$ specifies the contact order $c_{\rho}$ with respect to $\Delta_{\rho}$ if $\rho \in \sigma^{(1)}$, and $0$ otherwise.  

\begin{rema}\label{rem:contact-orders=lattice-points}
To any log smooth variety $X$ (where the underlying variety $\underline{X}$ is not necessarily smooth) one can associate a cone complex $\Sigma_X$ called the tropicalization of $X$. 
Then contact orders at a marked point are given precisely by the integral points of $\Sigma_X$. 
In the toric case $X = X_{\Sigma}$ as above, one checks that $\Sigma_X = \Sigma$ as cone complexes. In particular, the set of lattice points in $\Sigma$ are precisely the set of integral points of $\Sigma_X$. 

This canonical assignment of contact orders in the general situation is explained in \cite[\S 5.2]{ACMW17} and \cite[\S 2.3.8]{ACGS20}. 
\end{rema}

Henceforth we will identify a contact order with a lattice point in the support of the fan $|\Sigma|$ without any further mention.

\subsubsection{Campana curves on toric varieties}

Consider a log map $f \colon C \to X_{\Sigma}$ with the collection of contact orders $\varsigma = \{\mathbf{c}_k\}$. By intersection theory, this collection $\varsigma$ satisfies
\begin{equation}\label{eq:balancing-condition}
\sum_{\mathbf{c}_k \in \varsigma} \mathbf{c}_k = 0 
\end{equation}
which is called {\em the balancing condition}. Denote by $\mathbb{Z}\cdot \varsigma \subset N$ the sub-lattice generated by $\varsigma$.

\begin{theo}\label{thm:toric-uniruled-rc}
Let $\underline{X}_\Sigma$ be a smooth projective toric variety of dimension $d$ corresponding to the fan $\Sigma$ with a SNC divisor $\Delta_{\Sigma} \subset \underline{X}_\Sigma$.
Let $\varsigma$ be a collection of positive contact orders satisfying the balancing condition. Then 
\begin{enumerate}
 \item $X_{\Sigma}$ is separably $\varsigma$-uniruled. 
 \item If the sub-lattice $\mathbb{Z}\cdot \varsigma\subset N$ is of rank $d$, then $X_{\Sigma}$ is also rationally $\varsigma$-connected. 
 \item If $\mathbb{Z}\cdot \varsigma$ is of rank $d$, and $N/(\mathbb{Z}\cdot \varsigma)$ contains no $\ch\mathbf{k}$-torsion, then $X_{\Sigma}$ is separably rationally $\varsigma$-connected. 
\end{enumerate}
\end{theo}

\begin{proof}

Let $X_{\Sigma'}$ be another toric target of dimension $d$. We will argue that the above statements for $X_{\Sigma'}$ and for $X_{\Sigma}$ are equivalent. In particular, it suffices to consider the case that $\Sigma$ is the fan of a projective space, which then will be verified in Proposition \ref{prop:ps-uniruled-rc} below. 

Indeed, let $\pi : \mathcal U \to W, ev \colon \mathcal U \to X_{\Sigma}$
be a family of non-degenerate log maps over an irreducible $W$ with contact orders $\varsigma$ as in Definition \ref{def:uniruled-connected}. Note that $X_{\Sigma}$ and $X_{\Sigma'}$ share the same open dense tori $\mathbb{G}_m^d$ with the trivial log structure. 
For each geometric point $s \in W$, the non-degeneracy of the fiber $ev_s \colon U_s \to X_{\Sigma}$ implies that $U_s \setminus ev_s^{-1}(\mathbb{G}_{m}^d)$ consists precisely of markings. Hence $ev_s$ induces a non-degenerate log map $ev'_s \colon U_s \to X_{\Sigma'}$ of contact order $\varsigma$. Thus, there is an open dense $W' \subset W$ parametrizing non-degenerate log maps $ev' \colon \mathcal U \to X_{\Sigma'}$. 

If $X_{\Sigma}$ is separably $\varsigma$-uniruled, there is a family $W$ of genus zero non-degenerate stable log maps such that $\dim W = \dim X - 1$ and $ev$ is dominant and separable.  Then the corresponding family on $X_{\Sigma'}$ parametrized by $W'$ has the same properties, showing that $X_{\Sigma'}$ is $\varsigma$-uniruled.  Similarly, the properties of being rationally $\varsigma$-connected or separably rationally $\varsigma$-connected can be passed from $X_{\Sigma}$ to $X_{\Sigma'}$ using an appropriate choice of $W$ as in Definition \ref{def:uniruled-connected}.  Thus the desired statements for $X_{\Sigma'}$ follow from the analogous statements for $X_{\Sigma}$.
\end{proof}

\begin{rema}
In Theorem \ref{thm:toric-uniruled-rc} above, we restrict our attention to smooth toric varieties with SNC boundary because we have only defined ``$\varsigma$-uniruled'' and ``$\varsigma$-rationally connected'' in this setting.  We observe that the argument does not require that the fans $\Sigma$ and $\Sigma'$ are smooth; suitably interpreted, the claims are true for any (not necessarily smooth) complete fan $\Sigma$.  A key point is to interpret contact orders as lattice points of fans as in Remark \ref{rem:contact-orders=lattice-points}. 
\end{rema}

Since Theorem \ref{thm:toric-uniruled-rc} only imposes the the balancing condition on $\varsigma$, we can find suitable families of rational curves satisfying any multiplicity conditions.  For example:

\begin{coro}
\label{coro:toric_ratconn}
Suppose $(X, \Delta)$ is a projective klt Campana orbifold whose corresponding log scheme is $(X_{\Sigma},\Delta_{\Sigma})$ as in Theorem \ref{thm:toric-uniruled-rc}.  Then $(X,\Delta)$ is separably Campana rationally connected by good contact orders.
\end{coro}

\begin{proof}
Set $p = \mathrm{char}(\mathbf{k})$.  Let $m_{i}$ denote the multiplicity associated to the irreducible component $\Delta_{i}$ of $\Delta_{\Sigma}$.  Since $(X,\Delta)$ is klt, the constant $m = \sup_{i} \{m_{i}\}$ is finite.  Let $\alpha$ be a curve class such that $\Delta_{i} \cdot \alpha \geq 2m+2$ for every irreducible component $\Delta_{i}$.  Thus we can write $\Delta_{i} \cdot \alpha = t_{i,1} + t_{i,2}$ where (1) both $t_{i,1}$ and $t_{i,2}$ are at least $m$, and (2) neither $t_{i,1}$ or $t_{i,2}$ is divisible by $p$.

We define the contact order $\varsigma$ which has two markings for each torus invariant $\Delta_{i}$ and at these markings we assign either $t_{i,1}$ or $t_{i,2}$ to $\Delta_{i}$ and $0$ to every other torus-invariant divisor.  Since sum of the orders of $\varsigma$ along $\Delta_{i}$ agrees with $\Delta_{i} \cdot \alpha$ the contact order satisfies the balancing condition.  It is clear that $\mathbb{Z} \cdot \varsigma$ has full rank.  We claim that $N / \mathbb{Z} \cdot \varsigma$ has no $p$-torsion.  Indeed, suppose we fix a full-dimensional cone and consider the sublattice of $\mathbb{Z} \cdot \varsigma$ spanned by the vectors in $\varsigma$ proportional to these rays.  This already has no $p$-torsion, thus the same will be true for quotient by the entire sublattice $\mathbb{Z} \cdot \varsigma$.

By Theorem \ref{thm:toric-uniruled-rc} we see that $(X,\Delta)$ is separably rationally $\varsigma$-connected.  Since $\varsigma$ satisfies the Campana condition this finishes the proof. 
\end{proof}

\begin{proof}[Proof of Corollary~\ref{coro:toric}]
Weak approximation for a smooth projective toric variety over $k(B)$ is known by \cite[Theorem 4.3]{ColliotGille}. Thus our assertion follows from Theorem~\ref{thm:toric-uniruled-rc} combined with Theorem~\ref{theo:weakapproximation}.
\end{proof}

\subsection{Theorem \ref{thm:toric-uniruled-rc} for projective spaces}

\subsubsection{The set-up}\label{sss:projective-space-set-up}
For a positive integer $d$, consider the log scheme $X_d = X_{\Sigma}$ where the underlying $\underline{X_d} = \mathbb{P}^d$ is projective space. Let $\rho_0, \rho_1,\cdots, \rho_d$ be the ray generators of the rays in $\Sigma^{(1)}$. These satisfy $\sum_{i = 0}^d \rho_i = 0$. By an abuse of notation, we sometimes identify the generators $\rho_i$ with the rays that they span. 
The toric boundary is $\Delta_{X_d} = \sum_{j=0}^{d} H_j$ where $H_i$ is the hyperplane corresponding to $\rho_i$. 
The complement $X_d \setminus \Delta_{X_d} = \mathbb{G}_m^d$ is a rank $d$ torus.  We may choose homogeneous coordinate functions $[x_0:\cdots :x_d]$ of $X_d$ such that $H_{i} = (x_i = 0)$. 

\subsubsection{Parametrizing log curves in $X_d$}
We fix a rational curve $\underline{C} = \mathbb{P}^1$ with the homogeneous coordinate functions $[s:t]$ and set $b_{\infty} = [1:0] \in \underline{C}$. The complement $\underline{C} \setminus \{b_{\infty}\} = \mathbb{A}^1_{s/t} = \Spec \ \mathbf{k}[s/t]$ is an affine line. 

To parametrize rational log curves with the fixed underlying domain $\underline{C}$, we fix a positive integer $\beta$ to be the curve class, and a collection of non-zero lattice points $\varsigma_d = \{\mathbf{c}_k \in N\}_k $ where $\mathbf{c}_k$ specifies the contact order at the $k$-th marking $p_k$ as described in \S \ref{sss:toric-contact-orders}.
Note that as a lattice point in $N$, we have 
\[
\mathbf{c}_k = \sum_{i = 0}^{d} c_{k,i}\rho_i.
\]
The balancing condition  \eqref{eq:balancing-condition} is a consequence of the intersection-theoretic constraint
\[
\sum_{k=1}^{|\varsigma_d|} c_{k,i} = H_i.\beta
\]
for every $i$, which we now impose.

On $\underline{C}$ we fix a collection of markings 
\begin{equation}\label{eq:marking-corrdinates}
P := \{p_{k} = s_{k}/t_{k} \}_{k = 1}^{|\varsigma_d|} \subset \mathbb{A}^1_{s/t}.
\end{equation}
Suppose $f \colon C \to X_d$ is a non-degenerate rational log map with assigned contact orders $\varsigma_d$ at the markings $P$.
On the level of homogeneous coordinates, for each $i$ we have 
\begin{equation}\label{eq:log-curve-in-Pd}
f^*x_i = \lambda_i \cdot \prod_{k=1}^{|\varsigma_d|} (t\cdot s_{k}- s\cdot t_{k})^{c_{k,i}},
\end{equation}
for some $\lambda_i \in \mathbf{k}^{\times}$ since $f$ is non-degenerate. Inserting $b_{\infty} = [1:0]$ to \eqref{eq:log-curve-in-Pd}, we have 
\begin{equation}\label{eq:b-infty-position}
f^*x_i(b_{\infty}) = \lambda_i \cdot \prod_{k=1}^{|\varsigma_d|} (-t_{k})^{c_{k,i}}.
\end{equation}
Since $t_k \neq 0$, we observe that $(\lambda_i)_{i}$ hence $f$ is uniquely determined by the image $f(b_{\infty}) \in \mathbb{G}_m^d$.  Conversely, any such prescription determines a morphism $f: C \to X_{d}$ with the desired contact orders.

Denote by 
\begin{equation}
\label{equation:Ud}
U_{\varsigma} \subset (\mathbb{A}^1_{s/t})^{|\varsigma_d|} \times \mathbb{G}_m^d
\end{equation} 
the open subscheme parametrizing $|\varsigma_d|$ distinct points in $\mathbb{A}^1_{s/t}$ (representing the set of markings \eqref{eq:marking-corrdinates}) and a point $x_{\infty} \in \mathbb{G}_m^d$ (representing the image of $b_{\infty}$). Let $\pi \colon \mathcal C_{\varsigma} \to U_{\varsigma}$ 
be the family of log curves with markings specified by $U_{\varsigma}$. Note that $\underline{\mathcal C_\varsigma} \cong \underline{C}\times U_\varsigma$. We obtain a family of non-degenerate log curves
\[
f_\varsigma \colon \mathcal C_\varsigma \to X_d
\]
over $U_\varsigma$, fiberwise determined by \eqref{eq:log-curve-in-Pd}, such that $f_{\varsigma}(b_\infty) = x_{\infty}$. 

\begin{prop}\label{prop:ps-uniruled-rc}
Theorem \ref{thm:toric-uniruled-rc} holds for $X_d$. 
\end{prop}
\begin{proof}
The balancing condition and the parametrization \eqref{eq:log-curve-in-Pd} imply the existence of rational log curves in $X_d$ with contact order $\varsigma$. Theorem \ref{thm:toric-uniruled-rc} (1) for $X_d$ follows by using the $\mathbb G_m^d$-action.

\smallskip

Now consider the $2$-evaluation map
\[
ev^{(2)} := {f_\varsigma}\times_{U_\varsigma} {f_\varsigma} \colon \mathcal C_\varsigma\times_{U_\varsigma} \mathcal C_\varsigma \to X_d \times X_d.
\]
Assume that $\mathbb{Z}\cdot \varsigma$ is of rank $d$. 
To prove Theorem \ref{thm:toric-uniruled-rc} (2) for $X_d$, we will show that $ev^{(2)}$ is dominant.  It suffices to show that there is an $f$ defined as in \eqref{eq:log-curve-in-Pd} passing through a general pair of points $x, y \in \mathbb{G}_m^d$. Let $b_0 = [0:1] \in \underline{C}$. We will construct such $f$ satisfying
\begin{equation}
f(b_{\infty}) = x, \qquad  f(b_{0}) = y, \qquad b_0 \not\in P.
\end{equation}

Let $V^{\circ} \subset (\mathbb{A}^1_{s/t}\setminus \{b_0\})^{|\varsigma_d|}$ be the open subscheme parametrizing $|\varsigma_d|$ distinct points in $\mathbb{A}^1_{s/t}\setminus \{b_0\}$ as the set of markings \eqref{eq:marking-corrdinates}. Consider the locally closed subscheme $V_\varsigma:= V^{\circ}\times\{x\} \subset U_\varsigma$. By \eqref{eq:log-curve-in-Pd}, the restriction $(f_\varsigma)|_{\{b_0\}\times V_\varsigma} \colon \{b_0\}\times V_\varsigma \to X_d$ factors through $\mathbb{G}_m^d$, and is defined by 
\[
(f_\varsigma|_{\{b_0\}\times V_\varsigma})^*x_i = \lambda_i \cdot \prod_{k=1}^{|\varsigma_d|} s_{k}^{c_{k,i}}. 
\]
As $b_0, b_{\infty} \not \in P$, we may assume that $t_k = -1$ for all $k$ and $s_k \neq s_{k'}$ for $k\neq k'$.  
Hence for any $i$, $\lambda_i$ is the $i$-th homogeneous coordinate of $x$. 
Interpreting $(\mathbb{A}^1_{s/t}\setminus \{b_0\})^{|\varsigma_d|}$ as a torus, we see that the dimension of the image of  $(f_\varsigma)|_{\{b_0\}\times V_\varsigma}$ is given by the rank of $\mathbb{Z}\cdot \varsigma$, which is $d$. 
In particular, $(f_\varsigma)|_{\{b_0\}\times V_\varsigma} $ is dominant, implying that we can send $b_{0}$ to a designated general point $y$.

\smallskip

To prove Theorem \ref{thm:toric-uniruled-rc} (3) for $X_d$,  it suffices to show that   
\[
\operatorname{d} ev^{(2)} \colon T_{\mathcal C_\varsigma\times_{U_\varsigma} \mathcal C_\varsigma} \to (ev^{(2)})^*T_{X_d\times X_d}. 
\]
is surjective at the point $(b_0,u,b_{\infty}) \in \mathcal C_\varsigma\times_{U_\varsigma} \mathcal C_\varsigma$ for some general point $u \in U_\varsigma$. To compute 
$\operatorname{d} ev^{(2)}$, we choose the coordinate functions 
\[
\tilde{t}, \tilde{s}, \tilde{\lambda}_1,\cdots,\tilde{\lambda}_d, \tilde{s}_1,\cdots, \tilde{s}_{|\varsigma_d|}
\] 
around $(b_0,u,b_{\infty}) \in  \mathcal C_\varsigma\times_{U_\varsigma} \mathcal C_\varsigma$, and the coordinate functions
\[
\tilde{x}_1,\cdots, \tilde{x}_d, \tilde{y}_1,\cdots, \tilde{y}_d
\]
around the image $ev^{(2)}(b_0,u,b_{\infty}) \in X_d \times X_d$ as follows. 

Let $\tilde{t} = t/s$ (resp. $\tilde{s} = s/t$) be the coordinate around $b_\infty \in C_d$ (resp. $b_0 \in C_d$). Choose a general $u$ so that for any $p_k \in P$, we may assume $t_k = -1$, hence $\tilde{s}_k = s_k/t_k$ 
is the coordinate around $p_k \in \mathcal C_d$. In particular, we have $p_k \neq b_0$ for any $p_k \in P$, hence the image $ev^{(2)}(b_0,u,b_{\infty})$ avoids the boundary of $X_d$. 
Let $x_0, \cdots, x_d$ and $y_0, \cdots, y_d$ be the homogeneous coordinates of $X_d \times X_d$. We may assume that the coordinates around $ev^{(2)}(b_0,u,b_{\infty}) \in X_d \times X_d$  are given by $\tilde{x}_1 = x_1/x_0,\cdots, \tilde{x}_d = x_d/x_0$ and $\tilde{y}_1= y_1/y_0,\cdots, \tilde{y}_d = y_d/y_0$ and $x_0 = y_0 = 1$ at the point $ev^{(2)}(b_0,u,b_{\infty})$. Consequently by \eqref{eq:b-infty-position}, the coordinate around $f_\varsigma(b_{\infty}) \in \mathbb{G}_m^d$ is given by $\tilde{\lambda}_1 = \lambda_1/\lambda_0, \cdots, \tilde{\lambda}_d = \lambda_d/\lambda_0$ with $\lambda_0 = 1$ at $f_\varsigma(b_{\infty})$. 

Under the above choices of coordinates, using \eqref{eq:log-curve-in-Pd} we have 
\begin{equation}\label{eq:compute-ev-2}
(ev^{(2)})^*\tilde{x}_i = \tilde{\lambda}_i \cdot \prod_{k=1}^{|\varsigma_d|} (\tilde{s}_{k} +\tilde{s})^{c_{k,i}} \cdot \prod_{k=1}^{|\varsigma_d|} (\tilde{s}_{k} + \tilde{s})^{-c_{k,0}}, \qquad (ev^{(2)})^*\tilde{y}_i = \tilde{\lambda}_i \cdot \prod_{k=1}^{|\varsigma_d|} (\tilde{t}\cdot \tilde{s}_{k} +1)^{c_{k,i}}\cdot \prod_{k=1}^{|\varsigma_d|} (\tilde{t}\cdot \tilde{s}_{k} + 1)^{-c_{k,0}}
\end{equation}
around $(b_0,u,b_{\infty})$ for $i=1,\cdots,d$. 
Now the fiber $\operatorname{d} ev^{(2)}|_{(b_0,u,b_{\infty})}$ is given by evaluating the following $(2d)\times(d+|\varsigma_d|+2)$ Jacobian matrix
\[
\begin{bmatrix}
(\frac{\partial \tilde x_i}{\partial \tilde\lambda_j})_{i,j} &  (\frac{\partial \tilde x_i}{\partial \tilde s_k})_{i,k} & (\frac{\partial \tilde x_i}{\partial \tilde{s}})_{i} & (\frac{\partial \tilde x_i}{\partial \tilde{t}})_{i} \\
  (\frac{\partial \tilde y_i}{\partial \tilde\lambda_j})_{i,j} &  (\frac{\partial \tilde y_i}{\partial \tilde s_k})_{i,k} & (\frac{\partial  \tilde y_i}{\partial \tilde{s}})_{i} & (\frac{\partial \tilde y_i}{\partial \tilde{t}})_{i}
\end{bmatrix}
\]
at $(b_0,u,b_{\infty})$. 
Note that $|\varsigma_d| \geq (d+1)$ as $\mathbb{Z}\cdot \varsigma$ is of rank $d$, and $\varsigma_d$ satisfies the balancing condition. It suffices to verify the above matrix is of rank $2d$. 

A direct calculation shows that  
\[
 \left(\frac{\partial \tilde y_i}{\partial \tilde\lambda_j}\right)_{i,j}|_{(b_0,u,b_{\infty})} = I_{d\times d}
\]
where $ I_{d\times d}$ is the $d\times d$ identity matrix.
Next, we compute that 
\[
\frac{\partial \tilde x_i}{\partial \tilde s_k}|_{(b_0,u,b_{\infty})} = (c_{k,i} - c_{k,0})\frac{\lambda_i}{\tilde{s}_k} \cdot \prod_{\ell=1}^{|\varsigma_d|} \tilde{s}_{\ell}^{c_{\ell,i}-c_{\ell,0}} 
\]
As the factor $\frac{\lambda_i}{\tilde{s}_k} \cdot \prod_{\ell=1}^{|\varsigma_d|} \tilde{s}_{\ell}^{c_{\ell,i}-c_{\ell, 0}} \neq 0$, the rank of the matrix $(\frac{\partial \tilde x_i}{\partial \tilde s_k})_{i,k}|_{(b_0,u,b_{\infty})}$ is the same as the rank of the $d\times|\varsigma_d|$-matrix 
\begin{equation}\label{eq:contact-order-difference-matrix}
\Big(c_{k,i} - c_{k,0}\Big)_{i,k} \mod \ch\mathbf{k},
\end{equation} 
where $i$ runs through $\{1, 2, \cdots, d\}$ and $k$ runs through $\{1, 2, \cdots, |\varsigma_d|\}$. 

To compute the rank of \eqref{eq:contact-order-difference-matrix}, let $\tilde{N} \cong \mathbb{Z}^{d+1}$ be the lattice with generators $\{\tilde{\rho}_0, \cdots, \tilde{\rho}_d\}$. Consider the surjective morphism of lattices
\[
\varphi \colon \tilde{N} \longrightarrow N, \qquad \tilde{\rho}_i \mapsto  \rho_i,
\]
where $\rho_i$ is the lattice generator of the $i$-th ray of the fan $\Sigma^{(1)}$. 
Denote by $\tilde{\mathbf{c}}_k = \sum_{i} c_{k,i}\tilde{\rho}_i \in \tilde{N}$. Then observe that $\varphi(\tilde{\mathbf{c}}_k) = \mathbf{c}_k$ where we view $\mathbf{c}_k \in N$ as the lattice point. Let $\mathbb{Z}\cdot \tilde{\varsigma} \subset \tilde{N}$ be the sub-lattice generated by $\tilde{\varsigma} = \{\tilde{\mathbf{c}}_k\}_{k=1}^{|\varsigma|}$. 

Let $M$ and $\tilde{M}$ be the dual lattices of $N$ and $\tilde{N}$ respectively. Taking the dual of $\varphi$, we obtain an injection of lattices $\varphi^{\vee} \colon M \hookrightarrow \tilde{M}$. 
Let $\{\tilde{\rho}_i^{\vee}\}_{i=0}^{d}$ be the basis of $\tilde{M}$ dual to $\{\tilde{\rho}_i\}_{i=0}^{d}$. Then as a sub-lattice via $\varphi^{\vee}$, $M$ has a basis $\{\rho_i^{\vee}:= \tilde{\rho}_i^{\vee} - \tilde{\rho}_0^{\vee}\}_{i=1}^d$ dual to the basis $\{\rho_i\}_{i=1}^{d}$ of $N$. 
Observe that 
\[
c_{k,i} - c_{k,0} = (\tilde{\rho}_i^{\vee} - \tilde{\rho}_0^{\vee})(\mathbf{c}_k).
\]
Thus, the rank of \eqref{eq:contact-order-difference-matrix} is the dimension of 
\begin{equation}\label{eq:toric-dimension}
M|_{\mathbb{Z}\cdot \tilde{\varsigma}} = M|_{\mathbb{Z}\cdot {\varsigma}} \mod \ch\mathbf{k}.
\end{equation}
Here we view $M|_{*}$ as the space of linear functions defined on $*$. Finally, as $N/(\mathbb{Z}\cdot\varsigma)$ has no 
$\ch\mathbf{k}$-torsion, we observe that $M|_{\mathbb{Z}\cdot \tilde{\varsigma}}$ is of dimension $d$ as needed. 

This finishes the proof. 
\end{proof}

\begin{proof}[Proof of Theorem~\ref{theo:toric_intro}]
This follows from Theorem~\ref{thm:toric-uniruled-rc} and Corollary~\ref{coro:toric_ratconn}.
\end{proof}

\bibliographystyle{alpha}
\bibliography{CampanaCurve}

\end{document}